\documentclass[11pt,twoside,onecolumn]{article}

\usepackage{authblk} 

\usepackage{graphicx}
\graphicspath{ {Vis/} }
\usepackage[hang, small,labelfont=bf,up,textfont=it,up]{caption}
\usepackage{subcaption,dsfont}

\usepackage[usenames, dvipsnames]{color}

\usepackage{indentfirst} 

\usepackage[sc]{mathpazo} 
\usepackage[T1]{fontenc} 
\usepackage{microtype} 

\usepackage[english]{babel} 

\usepackage[hmarginratio=1:1,top=32mm,columnsep=20pt]{geometry} 
\usepackage{booktabs} 

\usepackage{lettrine} 

\usepackage{enumitem} 
\setlist[itemize]{noitemsep} 

\usepackage{titlesec} 

\usepackage{fancyhdr} 
\usepackage{hyperref} 

\usepackage{mathtools} 

\usepackage{amsmath,calc,amsthm}

\usepackage{algpseudocode}

\usepackage[numbers]{natbib}

\newcommand{\rev}{\textcolor{blue}}       
\newcommand\bff{\boldsymbol f}
\newcommand\bfF{\boldsymbol F}
\newcommand\bfa{\boldsymbol a}
\newcommand\bfe{\boldsymbol e}
\newcommand{\indi}[1]{{\rm I}_{\{#1\}}}

\theoremstyle{definition}

\theoremstyle{plain}
\newtheorem{proposition}{Proposition}

\usepackage{xcolor}

\usepackage{makecell}
\usepackage{scrextend} 

\usepackage{algorithm}
\usepackage{nomencl}
\makenomenclature
\usepackage{etoolbox}

\usepackage{multirow}

\makeatletter
\patchcmd{\thenomenclature}{\section*}{\section}{}{}
\makeatother

\begin{document}

\title{Dispatching Fire Trucks under Stochastic Driving Times} 

\author[1]{D. Usanov\thanks{Corresponding author. Tel.:+31(0)20 592 4168.\\
\textit{E-mail addresses}: usanov@cwi.nl (D. Usanov), p.m.van.de.ven@cwi.nl (P.M. van de Ven), r.d.van.der.mei@cwi.nl (R.D. van der Mei).}}
\author[1]{P.M. van de Ven}
\author[1, 2]{R.D. van der Mei}
\affil[1]{Centrum Wiskunde \& Informatica, Science Park 123,
1098 XG Amsterdam, The Netherlands}
\affil[2]{VU University Amsterdam, De Boelelaan 1105,
1081 HV Amsterdam, The Netherlands }

\date{} 


\maketitle
\newpage
\begin{abstract}
To accommodate a swift response to fires and other incidents, fire departments have stations spread throughout their coverage area, and typically dispatch the closest fire truck(s) available whenever a new incident arises. However, it is not obvious that the  policy of always dispatching the closest truck(s) minimizes the long-run fraction of late arrivals, since it may leave gaps in the coverage for future incidents. Although the research literature on dispatching of emergency vehicles is substantial, the setting with multiple trucks has received little attention. This is despite the fact that here careful dispatching is even more important, since the potential coverage gap is much larger compared to the single-truck case. Moreover, when dispatching multiple trucks, the uncertainty in the trucks' driving time plays an important role, in particular due to possible correlation in driving times of the trucks if their routes overlap.

In this paper we discuss optimal dispatching of fire trucks, based on a particular dispatching problem that arises at the Amsterdam Fire Department, where two fire trucks are send to the same incident location for a quick response. We formulate the dispatching problem as a Markov Decision Process, and numerically obtain the optimal dispatching decisions using policy iteration. We show that the fraction of late arrivals can be significantly reduced by deviating from current practice of dispatching the closest available trucks, with a relative improvement of on average about $20\%$, and over $50\%$ for certain instances. We also show that driving-time correlation has a non-negligible impact on decision making, and if ignored may lead to performance decrease of over $20\%$ in certain cases. As the optimal policy cannot be computed for problems of realistic size due to the computational complexity of the policy iteration algorithm, we propose a dispatching heuristic based on a queueing approximation for the state of the network. We show that the performance of this heuristic is close to the optimal policy, and requires significantly less computational effort.

\end{abstract}


\section{Introduction}

Due to the increased use of flammable synthetic materials in homes and offices, small fires may spread rapidly, potentially engulfing homes in a matter of minutes. In order to minimize property damage and save lives, many countries have strict laws that govern the fire departments' response time~\cite{legi}. To meet these requirements, fire departments operate a set of fire stations carefully located across their coverage area. When a new fire arises, one or more trucks are dispatched from the fire stations close to the fire in order to facilitate a quick response. However, sending closest trucks may lead to gaps in coverage for the duration of an incident which may have adverse effect on the response time to incidents that happen simultaneously. This is particularly true for large fires that require multiple trucks and take longer to put out. In this paper we study how to dispatch fire trucks in order to strike the right balance between responding quickly to the present fire, while maintaining good coverage for possible simultaneous incidents. 

To illustrate this tradeoff we consider the example of the Fire Department Amsterdam-Amstelland (FDAA), which operates 19 fire stations spread across the city of Amsterdam and surrounding areas. When a small fire occurs in the city center of Amsterdam that only requires a single fire truck to address, the FDAA nevertheless dispatches two trucks from different fire stations. These incidents are of the highest (of 3) priority level, and constitute about $70\%$ of all fires. When the first truck arrives at the fire, the second truck returns to its fire station. It is a policy FDAA uses for the city center where the streets are narrow, and in case there is a traffic jam, or an obstacle such as a garbage truck, the fire truck would not be able to overtake it but rather would have to go back and take an alternative route. Intuitively, the dispatcher would want to ensure that these two trucks are relatively close to the fire, but still sufficiently spaced out so that the remaining trucks retain good coverage. Moreover, we would want the trucks to approach the fire from different directions, so that when one truck gets stuck in traffic, the other can still get to the fire quickly. We refer to the latter phenomenon as {\em driving-time correlation}, and observe that this adds yet another layer of complexity to the optimal dispatching problem.

Although the problem of dispatching a single vehicle to incidents has been studied extensively in the literature on emergency services, to our knowledge very little work has been done on dispatching multiple vehicles, and we are the first to consider driving-time correlation in this context. Moreover, we are not aware of any studies into driving-time correlation in the transportation literature either. The current practice of the FDAA is to dispatch a truck each from the two fire stations closest to the incident. However, it is  unclear whether this leads to the fastest response (given the correlated driving times), and leaves the best coverage. Naturally, driving-time correlation also plays a role when considering incidents that require more than two trucks, but for ease of presentation we limit ourselves to the case with two trucks. While this problem is motivated by the situation of the FDAA, we believe other major cities with busy traffic use similar dispatching methods.

In order to study this problem, we model the city as a graph, where the vertices correspond to demand locations where incidents may occur, and an edge indicates that two locations can be reached directly. Fire stations are positioned at some of the vertices, and new fires arise at random times and locations. Similar to the current practice of the FDAA, we assume fires have to be addressed by sending two fire trucks, the first of which to arrive will engage the fire. \footnote{Note that we limit ourselves to the case of two trucks for simplicity, but we expect that our approach, heuristics and insights hold for larger fires that require more trucks.} The response time of a truck dispatched from a fire station to a fire is the sum of travel times over all edges traversed on the graph, and the travel time over each edge is some random variable. When two trucks dispatched to the same fire use the same edge they may incur the same travel times, capturing the driving time correlation.  Fires last for some random time, after which the trucks become idle again. In order to determine the optimal dispatching policy we model this system as a Markov decision process (MDP).

We first use policy iteration to numerically determine the optimal dispatching policy, and show that significant improvements can be made over the current practice of sending the two closest idle trucks. We also use this approach to demonstrate that it is important to take into account driving-time correlation in the model, since dispatch decision and performance metrics may be incorrect otherwise. For realistic-sized instances such as the coverage area of FDAA we cannot use policy iteration due to its computational complexity, and we develop novel heuristics instead.

Inspired by the results in~\cite{tiemessen2013dynamic}, we develop these heuristics using the idea of one-step improvement. This approach was developed in~\cite{norman1972heuristic,ott1992separable}, and has for instance been applied to call centers~\cite{bhulai2009dynamic}, control of traffic lights~\cite{haijema2008mdp}, routing in queueing networks~\cite{bhulai2003structure} and loss networks~\cite{hwang2000mdp}. To do this we first obtain an approximation for the fraction of late arrivals under the policy of sending the closest trucks, assuming that all fire stations are independent from each other. We then apply a single policy-iteration step to these results in order to obtain an improved policy. Although the independence assumption is very rough, we show that the resulting policy significantly outperforms closest-first. The computational complexity of this approach is much better than that of the full policy iteration algorithm needed to obtain the optimal dispatching policy, yet its performance is remarkably close to optimal. 

To summarize, in this paper we make the following contributions:
\begin{itemize}
\item[-] We develop the first model for dispatching multiple trucks in an emergency service network setting, possibly in the presence of correlated (stochastic) driving times;
\item[-] We show that the current fire department practice of sending the closest trucks is far from optimal, the optimality gap grows with the number of trucks in the system and can be as large as $50\%$ for certain problem instances;
\item[-] We show that taking into account driving time correlation has a significant impact on the response time and the optimal dispatch policy, and ignoring correlation when deriving a policy may lead to performance loss of more than $20\%$;
\item[-] To circumvent computational issues for obtaining the optimal dispatch policy, we propose a new heuristic based on 1-step policy improvement that has a small optimality gap, but only requires a fraction of its computational time.
\end{itemize}

In the next section we provide a review of the relevant literature, and in Section~\ref{sec:model} we give a description of the model studied in the paper, how we account for driving-time correlation, and formulate the MDP. In Section~\ref{sec:heuristics} we discuss one-step improvement policy and introduce our heuristics. In Section~\ref{sec:num} we numerically investigate the impact of correlation, compare the performance of the optimal policy, closest-first and the heuristics. Conclusions and suggestions for further research are made in Section~\ref{sec:conclusion}.

Throughout this paper we will denote vectors by boldfaced letters, e.g., $\boldsymbol x = (x_1,\dots,x_n)$, and by $|\boldsymbol x| = \sum_i |x_i|$ its 1-norm.

\section{Literature review}\label{sec:lit}

Operations research related to fire departments can be traced back to the RAND fire project, which ran from 1968 to 1975 and addressed a range of issues related to the New York City fire department. This includes for instance developing a simulation model for fire fighting services~\cite{carter1970simulation}, a square root law for fire fighting response times~\cite{kolesar1973square} and algorithms for relocations during major incidents~\cite{kolesar1974algorithm}. See~\cite{green2004anniversary} for an overview of this project and its research output. Since then the research literature on fire department operations has been limited in both scope and quantity, focussing mostly on facility location problems. The goal here is to determine the optimal location of the fire stations (see, e.g.,~\cite{marianov1992,kanoun2010,chevalier2012,degel2014}). 

To our knowledge the only papers that deal specifically with dispatching of fire trucks are~\cite{swersey1982markovian} and~\cite{ignall1982algorithm}, both originating from the RAND fire project. In~\cite{swersey1982markovian} the authors consider whether to dispatch one or two fire trucks to incidents of unknown severity, and show that the optimal policy has a threshold structure, where one only dispatches two trucks if there are sufficient trucks available. However, this paper ignores spatial components and does not determine {\em which} trucks to dispatch. The work closest to ours is perhaps~\cite{ignall1982algorithm}, where the authors propose an algorithm for how many (one or two) and which trucks to dispatch. \rev{The objective of the algorithm is to minimize response time to serious incidents, those requiring at least two ladder trucks.} The algorithm performs a grid search, where the first truck is picked for dispatching based on a certain loss approximation, assuming that only that truck is dispatched. Then, given the choice of the first truck, the second truck is decided on based on another loss function. Finally, the decision is made whether to send only the first truck or both of them based on the corresponding estimated costs. In contrast to our work,~\cite{ignall1982algorithm} relies on heuristic arguments for determining the future costs of current dispatching decisions, and ignores driving-time correlation. Moreover, the used loss functions do not seem to have an intuitive interpretation, and dispatching of the first truck is done independently of whether the second truck will be dispatched or not. In contrast, our approach is to jointly pick the two trucks to be dispatched such that the fraction of late arrivals is minimized, allowing to incorporate driving-time correlation. 

An area that is closely related to fire truck dispatching is that of dispatching ambulances to accidents and other emergencies. We will discuss the most relevant literature below, but emphasize that to our knowledge most of this work only considers dispatching a single vehicle to incidents, and does not take into account driving-time correlation. While results on the optimal dispatching of a single ambulance are not directly applicable to our setting, we now provide a brief discussion of some recent developments in this area. See for instance~\cite{ingolfsson2013ems,enayati2018ambulance, belanger2019} for a more complete overview of this field. 
\rev{In~\cite{andersson2007} a dispatching heuristic was proposed based on the notion of preparedness, measuring the ability of the system to respond quickly to future incidents. The heuristic suggests to dispatch an ambulance resulting in the smallest decrease in preparedness. The algorithm was further studied in~\cite{lee2011}. It was shown that the preparedness algorithm performs significantly worse than sending the closest ambulance in terms of average response time. The authors noted, however, that the poor performance of the preparedness algorithm is due to the fact that it ignores the current response time when making a dispatching decision. They introduced a modified version of the algorithm that balances between the decrease in preparedness and the response time to the current incident. In their experiments, the extended algorithm outperformed the closest-first dispatching policy.}

\rev{In~\cite{lim2011} the authors consider a setting with multiple incident priority levels, and compare a range of dispatching policies based on the closest-first policy. Modifications include possibilities to reroute busy ambulances to more urgent incidents and to reassign incidents to ambulances that become idle. The authors conclude that the relative performance of each policy depends on the parameters and available infrastructure.} In~\cite{jagtenberg2017dynamic} the authors formulate the problem of ambulance dispatching as an MDP, and then present a heuristic which is shown to perform close to optimal, and in certain cases outperforms closest-first. \rev{In~\cite{bandara2012} and~\cite{bandara2014} patient survivability is used as an objective for the problem with different incident priority levels. The authors formulate the problem as an MDP, and observe that dispatching closest vehicle is only optimal for the most urgent incidents. They also indicate that the optimal policy is most beneficial when the spacial distribution of incidents is unbalanced, which is the case in most real-life applications. Using the insights obtained from the optimal policy, the authors introduce a heuristic that outperforms the closest-first policy.} The authors of~\cite{mclay2013dispatching,mclay2013model} provide an MDP formulation of the ambulance dispatching problem under certain fairness constraints, and numerically compute the optimal policy for small instances. The problem of possibly sending two different types of emergency vehicles is considered in~\cite{sudtachat2014recommendations}, where the authors propose a heuristic for this purpose.

In addition to dispatching, substantial work in recent years has focused on relocation as well as joint dispatching and relocation of ambulances, in order to create better coverage. The relocation decisions imply proactive repositioning of idle vehicles within the region with the aim to reduce response time to future incidents. In~\cite{maxwell2010approximate,schmid2012,maxwell2013tuning,nasrollahzadeh2018real} the joint problem was addressed using approximate dynamic programming. In~\cite{enayati2018ambulance} the authors use stochastic programming to solve this problem, while ensuring that the workload due to relocations remains limited. The optimization method in~\cite{enayati2018real} is designed to make relocations that maximize coverage under personnel's workload limitations. Low computational costs of the approach allow to make decisions in real time, in contrast to the earlier methods described, which require offline computations.

\rev{As mentioned earlier, the research on ambulance dispatching is mostly focused on the setting where exactly one vehicle is required to serve an incident.} To understand why results for the single-vehicle case cannot easily be applied in our multiple-vehicle setting, consider the following. First, any dispatching action is a trade-off between a quick response, and ensuring that the remaining coverage is sufficient, should another incident arise while the first incident is still ongoing. Decomposing a multiple-vehicle dispatching formulation into a sequence of independent single-vehicle problems, one may not be able to carefully strike this balance, since every dispatching decision is made in a greedy way (assuming it is the only such decision). To illustrate this, consider the easier problem having to remove $k$ trucks: which set of $k$ trucks would result in the best coverage? It is easy to see that solving the problem sequentially would likely results in a substantially different solution with a worse coverage compared with solving the problem jointly for all trucks.

The second reason why algorithms for dispatching a single vehicle cannot be easily applied in our setting is due to the driving-time correlation. If applying single-vehicle policies for dispatching multiple trucks, one would be unable to take into account this correlation. As we shall show in this paper, driving-time correlation has a significant impact on the optimal dispatching policy, and ignoring it  substantially reduces performance.


\section{Model outline}\label{sec:model}
\begin{figure}[t]
	\centering
	\includegraphics[width=8cm]{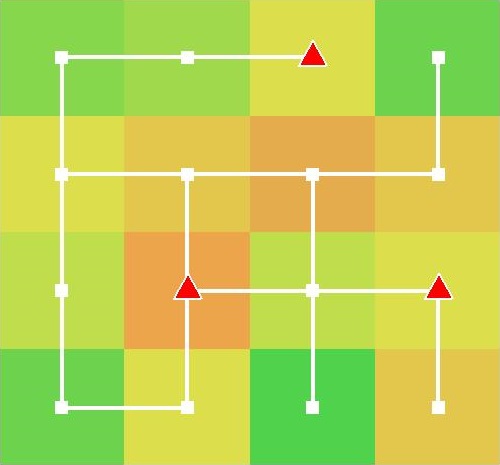}
	\caption{Graph representation of a region served by a fire department}
    \label{fig:graph}
\end{figure}

We consider a city represented by a connected, undirected graph $(\mathcal{J},E)$, see Figure \ref{fig:graph}. The set of vertices $\mathcal{J} = \{1,...,J\}$ represents the neighborhoods, or {\em demand locations}. Two vertices are connected if it is possible to travel directly between these two demand locations. A subset $\mathcal{I} \subseteq \mathcal{J}$ of demand locations contain a fire station (marked with triangles in Figure \ref{fig:graph}), and we denote  $I = |\mathcal{I}|$. Fire station $i \in \mathcal{I}$ houses $C_i$ fire trucks, and all fire trucks are assumed to be identical.

We assume that new fires arise at each demand location $j \in \mathcal{J}$ according to a Poisson process with rate $\lambda_j$, justified by the memorylessness of the time between new fires. Fire trucks can be either idle (waiting at a fire station) or busy (travelling or fighting a fire), and whenever a new fire starts, two idle fire trucks have to be dispatched. If a fire starts and fewer than two idle trucks are available, we request the missing truck(s) from a neighboring region. We assume the neighboring regions have ample capacity, so there are always trucks available. For tractability, we assume that when a truck is dispatched it remains busy for an exponential time with rate $\mu$, independent from the other truck dispatched and from the location of the fire and fire station. Independence from the location of the fire and fire station is a reasonable assumption as in practice the traveling time is negligible compared to the on-scene service time.
Note that the independence assumption allows us to consider any travel time distribution, although we shall focus mostly on Erlang distributed travel times, for ease of presentation. The assumption that both trucks have the same busy time distribution will result in an upper bound on the real-life busy fraction, since only the first truck to arrive will stay to resolve the incident. However, given the relatively low busy fraction for the fire truck application domain, we expect our model to be accurate. Returning trucks can be dispatched once they reach their station. Although an idealization, this assumption has negligible impact given the relatively low busy fraction of fire trucks seen in practice.

The state of the system can be represented by a vector $\bff = (f_1,\dots,f_I)$, where $f_i$ denotes the number of idle trucks at station $i$. Let $\boldsymbol a(\bff,j) = (a_1(\bff,j),...,a_I(\bff,j))$ represent the dispatch action taken if a new fire starts at a location $j$ when in state $\bff$. Here $0 \leq a_i(\bff,j)\leq f_i$ denotes the number of trucks dispatched from station $i \in \mathcal{I}$. Given that exactly two trucks are dispatched to every fire we have that $|\bfa(\bff,j)| \leq 2$, where the remaining $2-|\bfa|$ trucks are sent from neighboring regions.

We denote by $\bfF^{\bfa}(t)$ the state of the system at time $t$ under decision rule $\bfa$. Observe that, due to the exponentiality assumptions, the process $\{\bfF^{\bfa}(t)\}_{t\geq 0}$ is a continuous-time Markov process, with state space $\mathcal{S} = \{(f_1,...,f_I)|0 \leq f_i \leq C_i \ \forall i \in \mathcal{I}\}$, since each station $i$ can hold at most $C_i$ trucks. Let $\boldsymbol e_i$ denote a vector of length $I$ with $i$th element equal to 1, and all other elements equal to zero. The transition rates $q$ of this process are given by
\begin{align*}
q(\bff,\bff-\boldsymbol a(\bff, j)) &= \lambda_j,\quad j \in \mathcal{J}, \ \bff \in \mathcal{S},\\
q(\bff,\bff+\boldsymbol e_i) &= (C_i-f_i)\mu,\quad i \in \mathcal{I}, \ \bff \in \mathcal{S}.
\end{align*}
The first transition corresponds to trucks being dispatched upon the start of a new fire at location $j$, where the number of trucks at each location $i$ is reduced from $f_i$ to $f_i - a_i(\bff,j)$. These transitions occur at rate $\lambda_j$, the rate at which new incidents start at demand location $j$. The second transition corresponds to a truck returning to its fire station and becoming idle. This happens at rate $\mu$ for each individual truck not at its station, so the rate of trucks returning to station $i$ is equal to $(C_i-f_i)\mu$. This model resembles the hypercube model from~\citep{larson1974hypercube}. The hypercube model consists of a multiserver queueing model with distinguishable servers, corresponding to fire trucks in our setting. In~\cite{jarvis1981optimal} the authors numerically compute the optimal assignment policy of servers to requests in the hypercube model, and show that assigning the lowest-cost (closest in our setting) server is only optimal for small loads. The model is of relatively limited use in our setting, however, in that it cannot fully take into account the spatial component of our problem, and is only concerned with allocating a single server (dispatching a single truck). In~\cite{iannoni2007} a hypercube model was proposed used to analyze a system with particular dispatching policies including multiple dispatch and partial backup. This model was further embedded into a genetic algorithm in~\cite{iannoni2008} to optimize the service areas of ambulance bases.

\subsection{Traveling and response time}
When a truck is dispatched from fire station $i \in \mathcal{I}$ to demand location $j \in \mathcal{J}$, it travels along a shortest path on the graph, denoted by $s(i,j)$. Since we assume that the graph is connected, such a path always exists. In case multiple shortest paths exist, we select one at random. The travel time along edge $e \in E$ is denoted by $X_e \sim {\rm exp}(1)$, and follows an independent exponential random variable with unit mean. So the marginal traveling time of a fire truck dispatched from  $i$ to $j$ is given by $T_{i,j} = \sum_{e \in s(i,j)} X_e$, an Erlang distributed random variable with $|s(i,j)|$ phases of unit mean. 

For fire trucks dispatched from neighboring regions we assume a traveling time $T_0$ independent of the demand location of the fire, as typically those trucks are located relatively far and the driving time is dominated by the time it takes to reach the city in the first place. We assume that $T_0$ has an Erlang distribution with $2 \max_{i \in\mathcal{I},j \in \mathcal{J}} |s(i,j)|$ phases of unit mean. That is, the expected traveling time for a truck from a neighboring region is twice the maximum expected traveling time between any fire station-demand location pair on the graph, to reflect the fact that these trucks have to travel further.

The performance of a fire department is measured based on the response time to incidents, i.e., the time between the moment a fire reported and when the first truck arrives on scene. We consider two cases for computing the response time: \textit{uncorrelated} and \textit{correlated}. In the first case we use the simplifying assumption that the driving time on the same edge is independent between the two fire trucks. In the correlated case we assume that both trucks incur the same driving time realization for each shared edge. We now discuss each of these in more detail.

{\bf Uncorrelated driving times.} In order to model the fact that in the uncorrelated case the response times of the two trucks that are dispatched are completely independent, we introduce two independent copies of the driving time random variable over each edge. 
To do this we introduce an index $v=1,2$, which is used to distinguish between the two trucks that are dispatched, and is distinct from the index $i\in \mathcal{I}$ we use to index over all trucks. We denote  by $X_e^{(v)}$ the driving time of truck $v$ over edge $e \in E$, for $v=1,2$, and we assume that $X_e^{(1)}$ and $X_e^{(2)}$ are independent. We first treat the case where no trucks are sent from outside, and truck $v$ is dispatched from location $i_v$, $v=1,2$. In this case the total traveling time of the $v$-th truck to $j$ can be  written as $T_{i_v,j}^{(v)} = \sum_{e \in s(i_v,j)} X_e^{(v)}$, $v=1,2$. These $T_{i_v,j}^{(v)}$ are mutually independent because the $X_e^{(v)}$ are, even when $i_1=i_2$. The $T_{i_v,j}^{(v)}$ follow an Erlang distribution with $|s(i_v,j)|$ phases of unit mean.

In case one truck is dispatched from outside, we assume its traveling time is independent from the truck dispatched from inside the system; if two trucks are dispatched from outside their traveling times are assumed to be mutually independent. We denote by $T_0^{(1)}$ and $T_0^{(2)}$ two i.i.d.\ copies of the Erlang distributed random variable $T_0$.

Summarizing, in the uncorrelated case, given a dispatch decision $\bfa$ for a fire at location $j$, the response time can be expressed as
\begin{equation}\label{eqn:reponse_time_uncor}
R(\bfa,j) = \left\{
\begin{array}{ll}
\min\{T_{i,j}^{(1)},T_{i,j}^{(2)}\}		&	{\rm if~} a_i=2,\\
\min\{T_{i_1,j}^{(1)}, T_{i_2,j}^{(2)}\} 	&	{\rm if~} a_{i_1}=a_{i_2} = 1,\ i_1 \neq i_2,\\
\min\{T_{i,j}^{(1)},T_0^{(1)}\}			&	{\rm if~} a_i=1,\ |\bfa|=1,\\
\min\{T_0^{(1)},T_0^{(2)}\}				&	{\rm if~} |\bfa|=0.
\end{array}\right.
\end{equation}
The first two entries correspond to the case where two trucks are dispatched from inside the network with the first covering the case where both trucks are sent from the same location, and the second the case with different locations. Note that if the trucks are dispatched from the same station, they follow the same shortest path in a graph. This is a reasonable assumption as it is unlikely that in reality there are two independent shortest paths. Moreover, an alternative solution of sending the trucks via two different paths is hard to sell at the fire department as it is counterintuitive to the goal of getting to the incident as quick as possible. The third and fourth entry in~\eqref{eqn:reponse_time_uncor} correspond to the case where one and two trucks are dispatched from outside, respectively.

{\bf Correlated driving times.} In the correlated case the traveling times are no longer independent from each other, and we denote by $X_e$ the shared random traveling time over edge $e\in\mathcal{E}$ for both trucks. In contrast to the uncorrelated case, we need not distinguish between both trucks to compute the traveling time, and we denote $T_{i,j} = \sum_{e \in s(i,j)} X_e$ as the traveling time from $i$ to $j$ over $s(i,j)$, which is an Erlang distributed random variable with $s(i,j)$ phases. The traveling time of trucks dispatched from outside the network are still assumed to be independent from traveling times inside the network and from each other. Thus, in the correlated case the response time is given as follows:
\begin{equation}\label{eqn:reponse_time_cor}
R(\bfa,j) = \left\{
\begin{array}{ll}
T_{i,j}		&	{\rm if~} a_i=2,\\
\min\{T_{i_1,j}, T_{i_2,j}\} 	&	{\rm if~} a_{i_1}=a_{i_2} = 1,\ i_1 \neq i_2,\\
\min\{T_{i,j},T_0^{(1)}\}			&	{\rm if~} a_i=1,\ |\bfa|=1,\\
\min\{T_0^{(1)},T_0^{(2)}\}				&	{\rm if~} |\bfa|=0.
\end{array}\right.
\end{equation}
The entries correspond to the same decisions as in~\eqref{eqn:reponse_time_uncor} (respectively: two trucks from the same location, two trucks from different locations, one truck from outside the network, both trucks from outside the network). Note that in comparison to~\eqref{eqn:reponse_time_uncor}, the first entry no longer contains a minimum operator, since both trucks will have the same driving time realization as they are dispatched from the same location and there is correlation. The second entry is no longer necessarily a minimum  between two independent Erlang distributed random variables, as the routes of the two trucks may share one or more edges on the graph, for which they will see the same driving time realization.

Our approach described above for modeling driving-time correlation is certainly not the only possibility, and this work should be seen as the first attempt in taking this phenomenon into account when making dispatching decisions. For instance, note that we assume complete correlation between the driving time on each shared edge, whereas a smaller but still positive correlation coefficient may be more realistic. We briefly discuss this extension in Section~\ref{sec:conclusion}.

For each incident we are interested in whether the response time is within some time limit $t^*$, and we say a late arrival occurred otherwise. Our goal is to minimize the fraction of late arrivals. This is one of the most widely used performance metrics in emergency services, and is for instance used by the FDAA and the Dutch government to measure FDAA performance.

\subsection{MDP formulation}\label{sec:mdp}

We are interested in finding the dispatch decisions $\boldsymbol a(\bff,j)$ that minimize the fraction of late arrivals. In order to determine these we describe the system as an infinite-horizon average-cost Markov decision process (MDP). To do this we first uniformize our Markov process $\{\bfF^{\bfa}(t)\}_{t\geq 0}$ by adding the following dummy transitions: $q(\bff,\bff) = \mu \sum_{i \in \mathcal{I}} f_i$. This ensures that transitions out of any state happen at rate $\tau = \sum_{j\in \mathcal{J}}\lambda_j + \sum_{i \in \mathcal{I}}C_i$, without altering the dynamics of the network.

We are now in position to formulate our infinite-horizon average-cost MDP. Note that when a new fire starts and the network is in state $\bff$, we can make any of the following decisions $\bfa$:
$$
\mathcal{A}(\bff) = \{\bfa \in \mathds{N}_0^I \mid 0 \leq a_i \leq f_i,\ \min\{2,|\bff|\} \leq \sum_{i=1}^I a_i \leq 2\},
$$
i.e., we dispatch at most two trucks from inside the region, and we only dispatch outside trucks if fewer than two idle trucks are available. This description also states that we cannot dispatch more trucks from each station than available. Let $h^*(\bff)$ denote the relative cost incurred over an infinite time horizon when starting in state $\bff \in \mathcal{S}$, compared to paying the average cost $g^*$ every time unit. Since our process is unichain and has a finite state space and action space, we know from~\cite[Theorem 8.4.3]{puterman2014markov} that there exists an optimal deterministic policy that satisfies the Bellman equations,
\begin{align}
\nonumber h^*(\bff) \tau =& -g^* + \mu \sum_{i \in \mathcal{I}} (C_i - f_i) h^*(\bff +\bfe_i) + \mu \sum_{i \in \mathcal{I}} f_i h^*(\bff)\\
 &{}+ \sum_{j \in \mathcal{J}} \lambda_j \min_{\bfa \in \mathcal{A}(\bff)} \{\mathds{P}(R(\bfa,j)>t^*) + h^*(\bff-\bfa)\},\quad \bff \in \mathcal{S}. \label{eq:bellman}
\end{align}

The first summation on the right-hand side of~\eqref{eq:bellman} corresponds to fire trucks returning to their fire station, and the second to dummy transitions needed for uniformization. In neither case do we incur a cost or have to make a decision. The third summation corresponds to new fires that occur, in which case we have to make a dispatch decision $\bfa$, and incur some costs $\mathds{P}(R(\bfa,j)>t^*)$ equal to the probability of exceeding the response time threshold $t^*$, given the dispatch decision and location of the fire. The value function $g^*$ has an interpretation of the rate of late arrivals, that is, the average number of arrivals per time unit that were later than the time threshold $t^*$. To measure the performance of the dispatching policies we use the fraction of late arrivals, which is equal to $\frac{g^*}{\sum_{j \in \mathcal{J}} \lambda_j}$.

To compute the immediate costs $\mathds{P}(R(\bfa,j)>t^*)$, we must take a closer look at the distribution of the response time $R(\bfa,j)$, presented in~\eqref{eqn:reponse_time_uncor} and~\eqref{eqn:reponse_time_cor} for uncorrelated and correlated driving times, respectively. For uncorrelated driving times, in all four cases of~\eqref{eqn:reponse_time_uncor}, the response time is the minimum of two independent Erlang distributed random variables. The same holds for cases 3 and 4 of~\eqref{eqn:reponse_time_cor}, for correlated driving times.

The most challenging setting to compute is case 2 of~\eqref{eqn:reponse_time_cor}, where two trucks are dispatched from different locations under correlated driving times. This may be rewritten as the sum of an independent Erlang distributed random variable and the minimum of two others, i.e.,
\begin{equation}
R(\bfa,j) = \sum_{e \in s(i_1,j)\cap(i_2,j)} X_e + \min\{\sum_{e \in s(i_1,j)\setminus s(i_2,j)}  X_e, \sum_{e \in s(i_2,j)\setminus s(i_1,j)}  X_e\}, \ a_{i_1}=a_{i_2} = 1,\ i_1 \neq i_2.
\end{equation}
This kind of driving time correlation captures the fact that two fire trucks that take the same route may be delayed by the same incident or traffic, and encourages dispatching trucks over non-overlapping routes.

Thus, in order to compute the immediate costs $\mathds{P}(R(\bfa,j)>t^*)$, we require the following result.
\begin{proposition}\label{pro:Erlang}
Let $Y_0\sim Er(1, w_0)$, $Y_1\sim Er(1, w_1)$ and $Y_2\sim Er(1, w_2)$ be independent Erlang distributed random variables with phases of unit mean, $w_i>0$, $i=1,2,3$. Then
\begin{align*}
&P(\min\{Y_1, Y_2\} > t^*) = e^{-2t^*}\sum_{n=0}^{w_1-1}\sum_{m=0}^{w_2-1}\frac{{t^*}^{n+m}}{n!m!}\\
&and\\
&P(Y_0 + \min\{Y_1, Y_2\} > t^*)=\sum_{n=0}^{w_1-1}\sum_{m=0}^{w_2-1}\sum_{l=0}^{n+m}\frac{e^{-2t^*}{t^*}^l(-1)^{n+m-l}}{n!m!(w_0-1)!}\binom{n+m}{l}\int_{y_0 = 0}^{t^*}y_0^{n+m-l+w_0-1}e^{y_0}{\rm d}y_0\\
& + \sum_{n=0}^{w_0-1}\frac{{t^*}^n}{n!}e^{-t^*}.
\end{align*}
\end{proposition}

The proof of Proposition~\ref{pro:Erlang} can be found in Appendix~\ref{app:proof}.

\subsection{Closest-first dispatching}\label{sec:static}
The main benchmark throughout the paper is the current practice of FDAA, which is to always send the two closest (in terms of expected travel time) fire trucks, which we refer to as closest-first (CF) policy. We consider this as part of a larger class of static dispatching policies, where fire trucks are dispatched according to a fixed order per demand location. It can be represented by a list $\sigma_j(k)$, $j\in\mathcal{J}, k \in \{1,\dots,\sum_i C_i\}$, where $\sigma_j(k) \in \mathcal{I}$ represents the fire station from which to send the $k$th truck for an incident at location $j$. Let $\bfa^{CF}(\boldsymbol f,j)$ denote the action taken in state $\boldsymbol f$ given a new incident at location $j$, then

\begin{equation*}
\bfa^{CF}(\boldsymbol f,j) = \bfe_{\sigma_j(k_1)} + \bfe_{\sigma_j(k_2)},
\end{equation*}
where
$$
k_1 = \min\{k: f_{\sigma_j(k)}\ge 1\},	\quad	k_2 = \min\{k: f_{\sigma_j(k)}-\indi{k=k_1}\ge 1\},
$$
denote the number of the first and second truck dispatched, respectively. That is, truck $k_1$ is the closest fire truck to demand location $j$ that is currently present, and $k_2$ the second-closest. If $C_i=1$ for all $i$, then $\sigma_j$ reduces to a permutation over all fire stations. In case $k_i, i = 1, 2$ do not exist (because there are insufficient trucks available) we set $\sigma_j(k_i)=0$ and define $\bfe_0$ as the all-zero vector, to ensure trucks are sent from outside.

The long-term average costs under this CF policy can be obtained by limiting the Bellman equations~\eqref{eq:bellman} to only those actions $\bfa^{CF}(\boldsymbol f,j)$, i.e.,
\begin{align}
\nonumber h^{CF}(\boldsymbol f) \tau =& - g^{CF} + \mu \sum_{i \in \mathcal{I}} (C_i - f_i) h^{CF}(\bff +\bfe_i) + \mu \sum_{i \in \mathcal{I}} f_i h^{CF}(\bff)\\
 &{}+ \sum_{j \in \mathcal{J}} \lambda_j \left(\mathds{P}(R(\bfa^{CF}(\boldsymbol f,j),j)>t^*) + h^{CF}(\bff-\bfa^{CF}(\boldsymbol f,j))\right),\quad \bff \in \mathcal{S}. \label{eq:bellman_sd}
\end{align}
Here $g^{CF}$ and $h^{CF}(\boldsymbol f)$ denote the long-term average and relative costs under the CF policy, respectively. Thus~\eqref{eq:bellman_sd} is a system of $|S|$ linear equations, with $|S|+1$ unknowns $g^{CF}$ and $h^{CF}(\boldsymbol f)$, $f \in \mathcal{S}$. The costs can be obtained by fixing $h^{CF}(\boldsymbol f)$ for one state $\bff$, and solving the remaining system of equations.

\section{Dispatching heuristics}\label{sec:heuristics}

As we shall see from the experiments in Section~\ref{sec:num_opt}, the optimal dispatching policy significantly outperforms closest-first, both in the correlated and uncorrelated cases. However, it is well-known that solving the Bellman equations~\eqref{eq:bellman} can be computationally infeasible for large instances. In this section, we present two heuristics to approximate the optimal dispatching policy. 

\subsection{The OSI heuristic}

The first heuristic we consider is based on the idea of one-step improvement, and we refer to the policy obtained this way as to the one-step improvement (OSI) policy. This approach was developed in~\cite{norman1972heuristic,ott1992separable}, and the key idea is to first determine the (relative) costs $\tilde{h}(\boldsymbol y)$ for some sub-optimal policy, and then applying a single policy iteration step to find improved actions. That is, we replace the future costs $h^*(\boldsymbol y)$ in \eqref{eq:bellman} by some $\tilde{h}(\boldsymbol y)$. The maximizing action for this approximation of the Bellman equations can then be determined without iteration, significantly reducing the computational complexity compared to the full policy iteration algorithm. As pointed out in~\cite{norman1972heuristic,ott1992separable}, the first policy iteration step typically yields the biggest gains, so the result from one-step improvement is often close to optimal.

Here we use the CF policy to approximate the future optimal relative costs. We first compute the relative costs $h^{CF}(\boldsymbol f)$ from~\eqref{eq:bellman_sd}, and then substitute these into the right-hand side of the Bellman equations~\eqref{eq:bellman}. Ignoring the part that does not depend on the actions, the decision made by the OSI policy can be found as:
\begin{equation}\label{eq:osi}
\boldsymbol a^{OSI}(\boldsymbol f,j) \in \arg\min_{\boldsymbol a \in \textit{A}(\boldsymbol f)} \left(\mathds{P}(R(\bfa,j)>t^*) + h^{CF}(\bff-\bfa)\right),\quad \bff \in \mathcal{S}.
\end{equation}

\subsection{The OSIA heuristic}

To derive the OSI policy from \eqref{eq:osi}, we first need to solve the CF policy Bellman equations~\eqref{eq:bellman_sd} to determine the $h^{CF}(\bff)$. This is computationally expensive for large problem instances. In this section we present an algorithm that approximates the CF policy costs $h^{CF}(\bff)$, which can then in turn be used as a basis for the one-step improvement in~\eqref{eq:osi}. We will refer to the policy obtained using one step improvement with the CF policy cost approximation as the one-step improvement approximation (OSIA). This constitutes our second heuristic.

In order to approximate $h^{CF}(\bff)$, we assume that every fire station has exactly one truck. This assumption does not limit applicability of the algorithm, as we can always treat each truck as a separate station in the same location, and adjust the states and actions accordingly.

Let $J(\bff, t)$ denote the expected total cost under the CF policy during the time interval $[0,t]$ starting from state $\bff$. Then the relative cost $h^{CF}(\boldsymbol f)$ can be defined as
$$h^{CF}(\boldsymbol f) = \lim_{t \rightarrow \infty} \big( J(\boldsymbol f, t) - g^{CF}t \big),$$
where $g^{CF}$ denotes the cost per time unit under CF from~\eqref{eq:bellman_sd}.

Assume that after some time $T>0$ the system is in steady state, so the difference between the relative costs and the average costs is incurred in the interval $[0,T]$ only. In this case we can approximate $h^{CF}(\bff)$ as
\begin{align}
\nonumber h^{CF}(\boldsymbol f)  &= \lim_{t\rightarrow \infty} J(\bff,t)-g^{CF}t = J(\boldsymbol f, T) - g^{CF}T + \lim_{t\rightarrow \infty} (J(\bff,t) - J(\boldsymbol f, T))-(g^{CF}t - g^{CF}T)\\
&\approx J(\boldsymbol f, T) - g^{CF}T.  \label{eq:rel_sd_cost}
\end{align}

Substituting \eqref{eq:rel_sd_cost} into \eqref{eq:osi} we obtain the equations for the OSIA policy:
\begin{align*}
\boldsymbol a^{OSIA}(\boldsymbol f,j) &\in \arg\min_{\boldsymbol a \in \textit{A}(\boldsymbol f)} \mathds{P}(R(\bfa,j)>t^*) + J(\boldsymbol f-\bfa, T) - g^{CF}T\\
&= \arg\min_{\boldsymbol a \in \textit{A}(\boldsymbol f)} \mathds{P}(R(\bfa,j)>t^*) + J(\boldsymbol f-\bfa, T),\quad \bff \in \mathcal{S},\ j \in \mathcal{J}.
\end{align*}
Here we can omit the $g^{CF}T$ term because it appears for all actions.

So in order to derive the OSIA policy we need to estimate $J(\boldsymbol f, T)$, $\forall \boldsymbol f \in \textit{S}$, the total costs incurred in the interval $[0,T]$, starting from state $\bff$. Following an idea from~\cite{tiemessen2013dynamic}, we decompose the network into individual $M/M/1/1$ queues associated with individual fire stations. By doing this, we essentially decouple the network into individual fire stations, for each we can now compute an approximation for the probability of the corresponding fire truck to be busy (busy probability). These we combine to obtain an approximation for $J(\boldsymbol f, T)$.

Let us first consider a fire station $i$ in isolation, and compute its busy probability. Denote by $D_i$ the given demand arrival rate for the truck at station $i$. Recall that the steady-state busy probability of an $M/M/1/1$ queue with load $\rho_i$ is given by $B(\rho_i) = \rho_i/(1+\rho_i)$, and thus the steady-state rate of rejected requests is $D_i B(\rho_i)$. Denote by $N(\rho_i, f_i, t)$ the expected number of rejected requests in the $M/M/1/1$ queue during $[0, t]$ starting with $f_i$ trucks at time $0$. Finally, let $\Delta(\rho_i, f_i)$ be the difference in rejected requests between starting from steady state and starting from $f_i$: $\Delta(\rho_i, f_i) = \lim_{t \rightarrow \infty} \big(N(\rho_i, f_i, t) - D_i t B(\rho_i) \big)$.

Assuming as above that the system is in steady state after time $T$, we have that 
\begin{equation}\label{eqn:diff}
\Delta(\rho_i, f_i) \approx \big(N(\rho_i, f_i, T) - D_i T B(\rho_i) \big).
\end{equation}
The busy probability $p_i$ can be obtained by dividing the expected total number of rejections $N(\rho_i, f_i, T)$ by the expected number of arrivals $D_i T$. Observe that in our case $\rho_i  = D_i/\mu$, since each request will occupy the server (i.e., fire truck) for an expected duration $\mu^{-1}$. Using the identity in~\eqref{eqn:diff} and bounding between 0 and 1 to obtain a probability (since we are using approximations), we get
\begin{equation}\label{eq:busy_prob}
p_i = \frac{N(\rho_i, f_i, T)}{D_i T} = \max \big\{0, \min \big\{1, B(D_i/\mu) + \frac{\Delta(D_i/\mu, f_i)}{D_i T} \big\} \big\}.
\end{equation}

Observe that in order to evaluate~\eqref{eq:busy_prob} we need to approximate $\Delta(D_i/\mu, f_i)$, the difference in total number of rejected calls between steady-state and starting from state $f_i$. To do this, we formulate the queue representing station $i$ as an average-cost MDP, where the state is the number of idle trucks at the fire station. Transitions happen when either a request for a truck arrives or an idle truck returns from an incident. If there is an idle truck, it is always dispatched. The cost for a rejection is 1, and 0 for an accepted job. This results in the following system of 2 Bellman equations and a normalizing equation:
\begin{align}
  &h_0 = \frac{D_i}{D_i + \mu} - \frac{D_iB(\frac{D_i}{\mu})}{D_i + \mu} + \frac{D_i}{D_i + \mu}h_0 + \frac{\mu}{D_i + \mu}h_1, \label{eq:mm11_1}\\
  &h_{1} = -\frac{D_iB(\frac{D_i}{\mu})}{D_i} + h_{0}, \label{eq:mm11_2}\\
  &\frac{1}{1+\frac{D_i}{\mu}}h_0 + \frac{\frac{D_i}{\mu}}{1 + \frac{D_i}{\mu}}h_1=0. \label{eq:mm11_3}
\end{align}
Solving \eqref{eq:mm11_1}-\eqref{eq:mm11_3}, we obtain $\boldsymbol h = (h_0, h_1)$, the relative costs starting from state $f_i=0$ or $f_i=1$, respectively. We use $\Delta(D_i/\mu, f_i)=h_{f_i}$, and compute $p_i$ using \eqref{eq:busy_prob}. 

Having determined the busy probability $p_i$ for a given arrival rate $D_i$, our next step is to update the values of $D_i$ using the busy probabilities obtained. Here we again consider all fire stations jointly. According to the CF policy, the closest two idle trucks are dispatched to an incident. Recall that the lists $\sigma_j(k)$, $j\in\mathcal{J}, k \in \{1,\dots,I\}$ represent the dispatching order corresponding to the CF policy. So as each station has exactly one truck, $\sigma^{-1}_j(i)$ denotes the position held by station $i$ in the dispatching order of demand location $j$. For instance, $\sigma^{-1}_j(i) = 1$ means that station $i$ is the closest to demand location $j$.

Let $p_0$ correspond to the probability of an outside truck being unavailable, and set $p_0=0$.  
After $p_i$ is computed for each station $i$ according to~\eqref{eq:busy_prob}, we calculate the probability $p^j_{\{i_1,i_2\}}$ of a newly arrived incident at demand location $j$ requests trucks at $i_1$ and $i_2$. Note that a single incident can generate requests at multiple pairs of fire stations, since some of them might be occupied. By conditioning on the availability of the fire trucks we obtain:\\
For $j=1,..,J$, $i_1 = 2,...,I$, $i_2 = 1,..,(i_1-1)$ (both trucks are from inside):
\begin{equation}\label{eq:couple_prob}
p^j_{\{i_1,i_2\}} =
\begin{cases}
	1 & \mbox{, if } \sigma_j(i_1) = 1,\ \sigma_j(i_2) = 2, \\
	&\mbox{  or } \sigma_j(i_1) = 2,\ \sigma_j(i_2) = 1,\\
	
	\prod_{i \neq i_1, i_2,\ \sigma_j(i)<max\{\sigma_j(i_1), \sigma_j(i_2)\}}p_i & \mbox{, otherwise.}
\end{cases}	
\end{equation}
For $j=1,..,J$, $i_1 = 1,...,I$, $i_2 = 0$ (one truck is from outside):
\begin{equation}\label{eq:couple_prob_0_i}
p^j_{\{i_1,i_2\}} = \prod_{i \neq i_1}p_i.
\end{equation}
For $j=1,..,J$, $i_1 = 0$, $i_2 = 0$ (both trucks are from outside):
\begin{equation}\label{eq:couple_prob_0_0}
p^j_{\{i_1,i_2\}} = \prod_{i = 1}^{I}p_i.
\end{equation}

The probability $p^j_{\{i_1,i_2\}}$ is equal to 1 if trucks at $i_1$ and $i_2$ are the closest to $j$. Otherwise, it is a product of the busy probabilities of those trucks that are closer than either $i_1$ or $i_2$. Trucks from inside of the region are always closer than those from outside.

Denote $D_{\{i_1,i_2\}}$ the demand arriving for trucks from stations $i_1$ and $i_2$. Given the probabilities $p^j_{\{i_1,i_2\}}$, we compute $D_{\{i_1,i_2\}}$ for $i_1=2,...,I, \ i_2=1,...,(i_1-1)$:
\begin{equation}\label{eq:couple_demand}
D_{\{i_1,i_2\}} = \sum_{j \in \mathcal{J}}\lambda_j p^j_{\{i_1,i_2\}}.
\end{equation}
Finally, by summing over all pairs $\{i,k\}$, $k \neq i$, we can obtain the arrival rate of incidents at station $i$ as $D_i = \sum_{k\neq i}D_{\{i,k\}}$.

Let $C^j_{i_1i_2}$ indicate the expected penalty related to sending trucks $i_1$ and $i_2$ to location $j$. It is equivalent to the cost 
$\mathds{P}(R(\bfa,j)>t^*)$ where the action $\boldsymbol a$ corresponds to sending the trucks from stations $i_1$ and $i_2$ to location $j$, given that those are idle. Costs computation is discussed earlier in Section \ref{sec:mdp}. We now summarize the algorithm that approximates $J(\boldsymbol f, T)$ for a given state $\boldsymbol f \in \mathcal{S}$ in pseudocode Algorithm \ref{alg:approx}.

\begin{algorithm}[t]
\caption{CF cost approximation}\label{alg:approx}
\begin{algorithmic}
\State \textbf{Initialization}
\State $p^j_{\{i_1,i_2\}} =
\begin{cases}
	1, & \mbox{if } \sigma_j(i_1) = 1,\ \sigma_j(i_2) = 2 \mbox{ or } \sigma_j(i_1) = 2,\ \sigma_j(i_2) = 1\\
	
	0, & \mbox{otherwise}
\end{cases}	$
\State $D_{\{i_1,i_2\}} = \sum_{j \in \mathcal{J}}\lambda_j p^j_{\{i_1,i_2\}}\ \ \ \forall i_1, i_2 \in \{0,1,...,I\}$
\State $D_i = \sum_{k\neq i}D_{\{i,k\}} \ \ \ \forall i \in \mathcal{I}$
\While{true}
\State{Compute $\Delta(D_i, \mu, f_i)=h_{f_i}$ using \eqref{eq:mm11_1}-\eqref{eq:mm11_3}}
\State{Compute $p_i$ using \eqref{eq:busy_prob}}
\State{Compute $p^j_{\{i_1,i_2\}}$ using \eqref{eq:couple_prob}-\eqref{eq:couple_prob_0_i}}
\State{Compute $D_{\{i_1,i_2\}}$ using \eqref{eq:couple_demand}}
\State $\hat{D_i} = \sum_{k\neq i}D_{\{i,k\}} \ \ \ \forall i \in \mathcal{I}$
\If{$|D_i-\hat{D_i}|/D_i < \epsilon \ \ \ \forall i \in \mathcal{I}$}
\State $D_i = \hat{D_i} \ \ \  \forall i \in \mathcal{I}$
\State{break}
\EndIf
\State $D_i = \hat{D_i} \ \ \  \forall i \in \mathcal{I}$
\EndWhile

\State $J(\boldsymbol f, T) = T\sum_{j\in\mathcal{J}}\lambda_j \sum_{i_1 = 0}^{I}\sum_{i_2=0}^{\max\{0,i_1-1\}} p^j_{i_1i_2} C^j_{i_1i_2} (1 - p_{i_1})(1 - p_{i_2})$
\end{algorithmic}
\end{algorithm}

\section{Numerical results}\label{sec:num}

We now present the results of our numerical experiments. In Section~\ref{sec:setup} we describe the setup of our numerical experiments. The results are separated into two parts: in Section~\ref{sec:num_opt} we compare the CF and OPT policies, and use this to understand how much improvement over CF can be obtained, and what is the impact of driving-time correlation on the policies and their performance. In Section~\ref{sec:heur} we then evaluate the performance of our heuristics OSI and OSIA relative to CF and OPT, both in terms of fraction of late arrivals and computational time.

\subsection{Setup of the numerical experiments}\label{sec:setup}

All experiments were run in MATLAB R2017b on a computer with an Intel Core i5-5250U 1.6 GHz processor, 8 GB RAM, running Linux Fedora 26. In order to evaluate the performance of a policy for a given network and set of parameters, we numerically solve the Bellman equations~\eqref{eq:bellman} for OPT policy and the restricted Bellman equations~\eqref{eq:bellman_sd} for CF policy. This way we obtain $g^{OPT}$ and $g^{CF}$, the long-term expected number of late arrivals per time unit for OPT and CF, respectively. The dispatching order $\sigma_j(k)$ for CF is determined by ordering for each demand location $j$ the fire stations $k$ based on the length of their shortest path to $j$. Ties are broken arbitrarily. 

In order to compute the performance of OSI we first determine the relative costs for closest first $h^{CF}(\boldsymbol f)$ from~\eqref{eq:bellman_sd}, and substitute these into~\eqref{eq:osi} to determine the actions $\boldsymbol a^{OSI}$. These are then substituted into the Bellman equations~\eqref{eq:bellman}, which we solve numerically to obtain the rate of late arrivals for OSI $g^{OSI}$. For OSIA we repeat this procedure, except that instead of computing the exact relative costs for closest first $h^{CF}(\boldsymbol f)$, we compute $J(\boldsymbol f,T)$ from Algorithm~\ref{alg:approx} and use the approximation for $h^{CF}(\boldsymbol f)$ from~\eqref{eq:rel_sd_cost}. This way we obtain $g^{OSIA}$, the rate of late arrivals under OSIA. In order to compute the fraction of late arrivals (FLAR) for any of these policies, we divide the long-term expected number of late arrivals per time unit $g$ by the total arrival rate, i.e., $g/\sum_{j \in \mathcal{J}} \lambda_j$.

\begin{figure}
\begin{subfigure}{.3\textwidth}
\centering
\includegraphics[width=\textwidth]{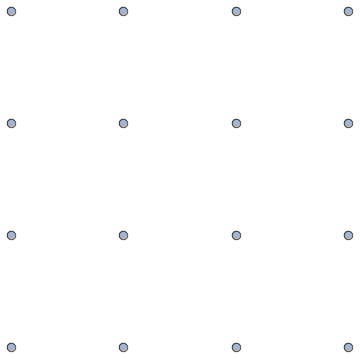}
\caption{}
\label{fig:grid_graph_1}
\end{subfigure}\hfill
\begin{subfigure}{.3\textwidth}
\centering
\includegraphics[width=\textwidth]{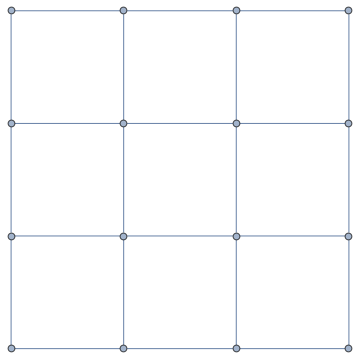}
\caption{}
\label{fig:grid_graph_2}
\end{subfigure}\hfill
\begin{subfigure}{.3\textwidth}
\centering
\includegraphics[width=\textwidth]{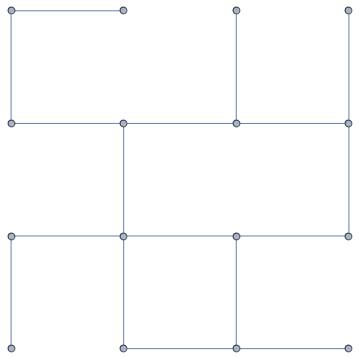}
\caption{}
\label{fig:grid_graph_3}
\end{subfigure}
\caption{Random graph construction}
\end{figure}

For our experiments we randomly generate grid-like graphs, as outlined below. For some parameter $d \in \mathbb{N}$, we generate a grid of $d \times d$ vertices (see Figure~\ref{fig:grid_graph_1}), placed at unit distance. We then connect each pair of vertices within unit distance from each other, so a vertex away from the boundary is connected to its four immediate neighbors (see Figure~\ref{fig:grid_graph_2}) and we obtain a graph with $|\mathcal{J}| = d^2$ nodes and $|E| = 2d(d-1)$ edges. We then remove edges uniformly at random until the number of removed edges is below $2d(d-1) s$ (see Figure~\ref{fig:grid_graph_3}), where $s \in (0,1)$ is some desired level of sparseness. The $s$ parameter is drawn from a uniform distribution $\mathcal{U}(0.4, 1)$. While removing the edges, we check if the graph remains connected. In case the graph becomes disconnected, a new random edge is selected for removal. If after a certain number of attempts no edge is found that can be removed without disconnecting the graph, the procedure stops, and the obtained graph is used.

In our experiments we assume each station has exactly one truck. This does not affect methodology, but makes it easier to visualise and understand the difference in actions taken by different policies. We allocate stations (or trucks) to vertices sequentially in a randomized manner. Each of the $I$ trucks is positioned on a vertex not yet occupied by other trucks uniformly at random.


\subsection{Comparison of closest-first and optimal dispatching}\label{sec:num_opt}

In this section we are interested in studying OPT and its performance relative to the CF heuristic. Recall that OPT is computed from the Bellman equations~\eqref{eq:bellman} through policy iteration, and it is here that we run into the so-called curse of dimensionality, which states that the state space and action space of the MDP become too big to solve in an efficient manner. Specifically, our action space grows as $\mathcal{O}(I^2)$ since each action consists of sending two trucks. The state space grows as $\mathcal{O}(2^I \times d^2)$, since there are $2^I$ possible combinations of available trucks, and the next fire can occur on any of the $J=d^2$ demand locations. Although the complexity of each step of policy iteration is polynomial in the size of the state space and action space, there is no universal polynomial bound on the complexity of the algorithm, due to the uncertainty in the number of steps required~\cite{littman1995complexity}. In practical terms, this means that we can only compute the optimal policy for instances of small-to-moderate size. In Section~\ref{sec:heur} we restrict ourselves to suboptimal policies, and consider instances of real-life size (in the case of FDAA   there are roughly $I=13$ trucks and $J \approx 400$ demand locations). Due to the relatively low load seen in the FDAA practice ($\rho = 0.02$) and used in our experiments, the number of incidents that requires trucks from outside is negligible.

\begin{figure}[t]
\begin{subfigure}{.5\textwidth}
\centering
\includegraphics[width=\textwidth]{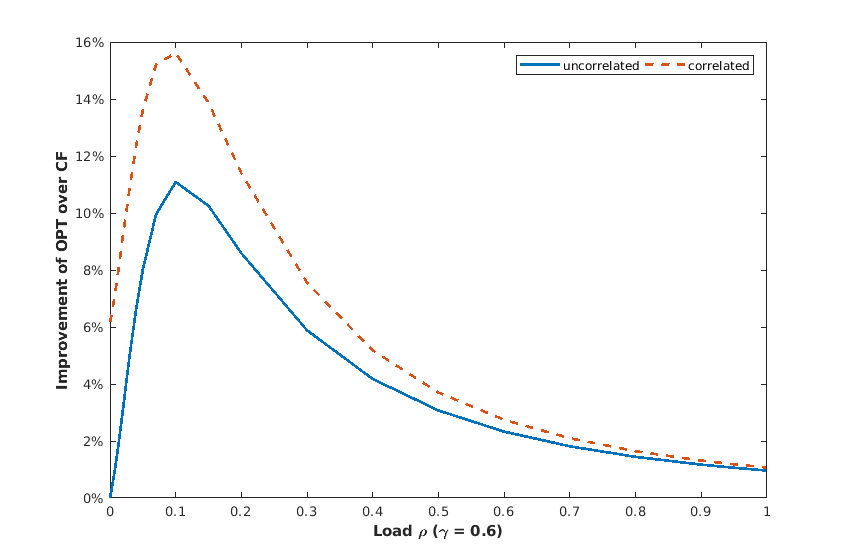}
\caption{graph 1}
\label{fig:varying_rho_1}
\end{subfigure}\hfill
\begin{subfigure}{.5\textwidth}
\centering
\includegraphics[width=\textwidth]{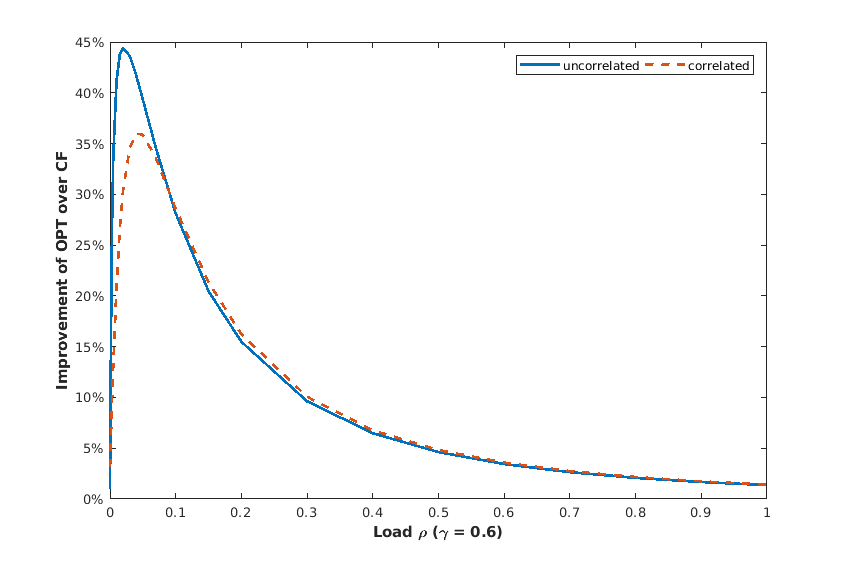}
\caption{graph 2}
\label{fig:varying_rho_2}
\end{subfigure}\hfill
\begin{subfigure}{.5\textwidth}
\centering
\includegraphics[width=\textwidth]{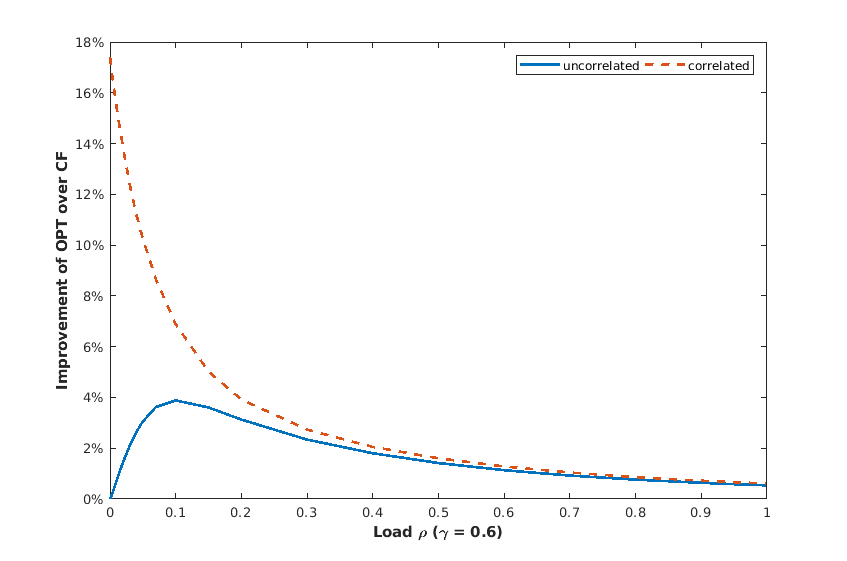}
\caption{graph 3}
\label{fig:varying_rho_3}
\end{subfigure}
\begin{subfigure}{.5\textwidth}
\centering
\includegraphics[width=\textwidth]{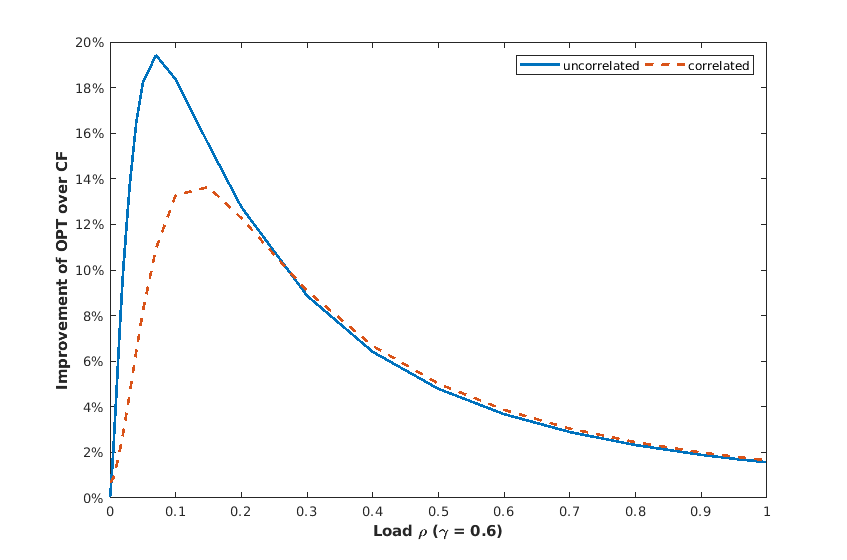}
\caption{graph 4}
\label{fig:varying_rho_4}
\end{subfigure}
\caption{$\delta^{OPT}$ as a function of $\rho$ for 4 random graphs ($I = 5$, $d = 7$, $\gamma = 0.6$)}
\label{fig:varying_rho}
\end{figure}

{\bf Relative improvement of the optimal policy over closest-first.} We are interested in assessing the current practice of dispatching the two closest trucks, and to see whether there is any room for improvement (i.e., reducing the fraction of late arrivals) by dispatching in a smarter way. To do this, we consider the relative improvement of OPT over CF, which is computed as
$$\delta^{OPT} = \frac{g^{CF}-g^{OPT}}{g^{CF}}\times 100\%.$$
In Figure~\ref{fig:varying_rho} we plot the relative improvement against the load of the system $\rho = \frac{\sum_{j \in \mathcal{J}}\lambda_j}{I\mu}$, which represents the amount of work per truck arriving each time unit. We do this for four different randomly generated graphs, and show the improvement both in uncorrelated and correlated cases. We define the time threshold for late arrivals as $t^* = \gamma \max_{i \in \mathcal{I}, j \in \mathcal{J}}{|s(i,j)|}$, to ensure that it scales with the graph size, and set $\gamma=0.6$.

We see that in both cases the relative improvement depends on the graph, and ranges from $0\%$-$50\%$, depending on the load and on this graph. This is significant, and suggests that in the right circumstances, significant gains can be found by dispatching in a clever way. In the uncorrelated case the relative difference is small when $\rho$ is small or large. This is because if the load is close to 0, the system is almost always in the state with all the trucks being idle, and when $\rho$ is close to 1, there is no room for improvement independent of whether there is correlation or not, because the system is almost always in the state with no idle trucks.

When correlation is introduced however, we see from Figure \ref{fig:varying_rho} that sending two closest trucks does not necessarily minimize response time, even for small loads. Hence, in this case the OPT policy may improve upon the CF policy even for very small values of $\rho$, as illustrated in Figures~\ref{fig:varying_rho_3} and~\ref{fig:varying_rho_4}. 
However, we see in all cases in Figure~\ref{fig:varying_rho} that as $\rho$ grows, the improvement curve foro uncorrelated driving times converges to the one corresponding to uncorrelated case.
 
%

\begin{figure}[t]
\begin{subfigure}{.5\textwidth}
\centering
\includegraphics[width=\textwidth]{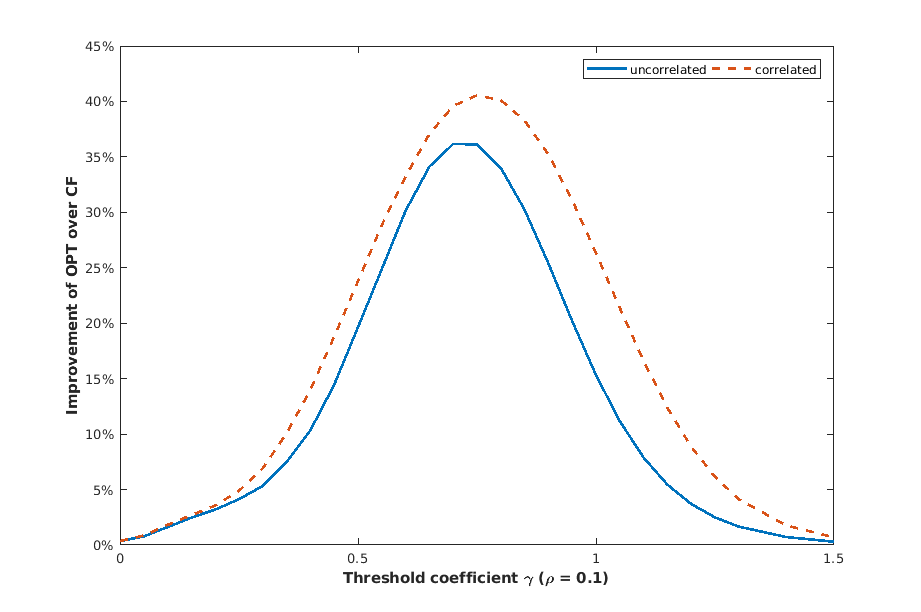}
\caption{graph 1}
\label{fig:varying_thr_1}
\end{subfigure}\hfill
\begin{subfigure}{.5\textwidth}
\centering
\includegraphics[width=\textwidth]{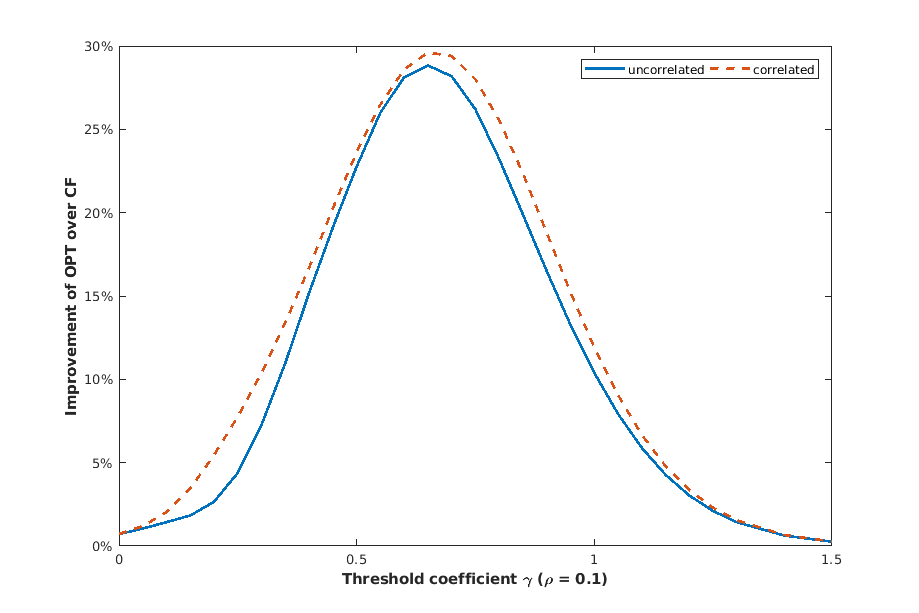}
\caption{graph 2}
\label{fig:varying_thr_2}
\end{subfigure}\hfill
\begin{subfigure}{.5\textwidth}
\centering
\includegraphics[width=\textwidth]{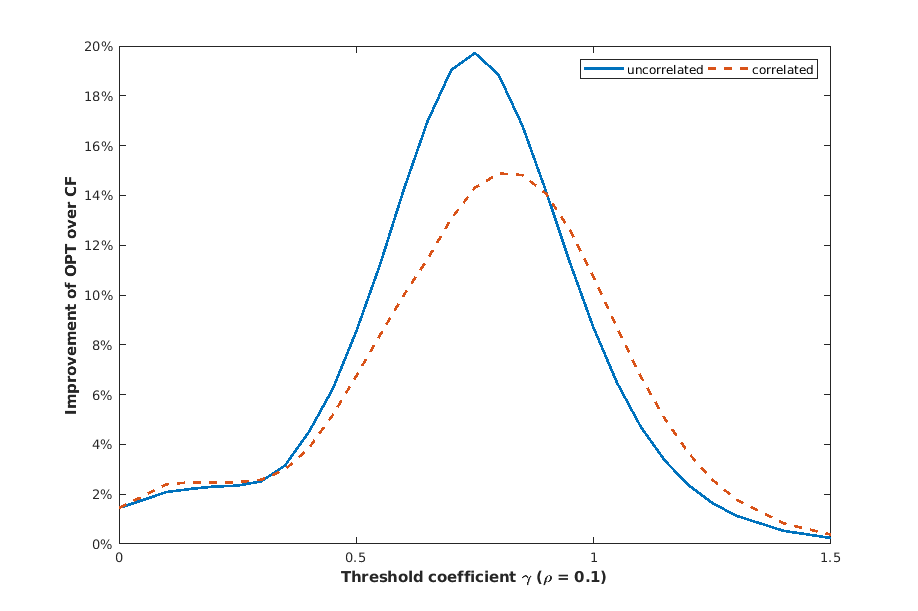}
\caption{graph 3}
\label{fig:varying_thr_3}
\end{subfigure}
\begin{subfigure}{.5\textwidth}
\centering
\includegraphics[width=\textwidth]{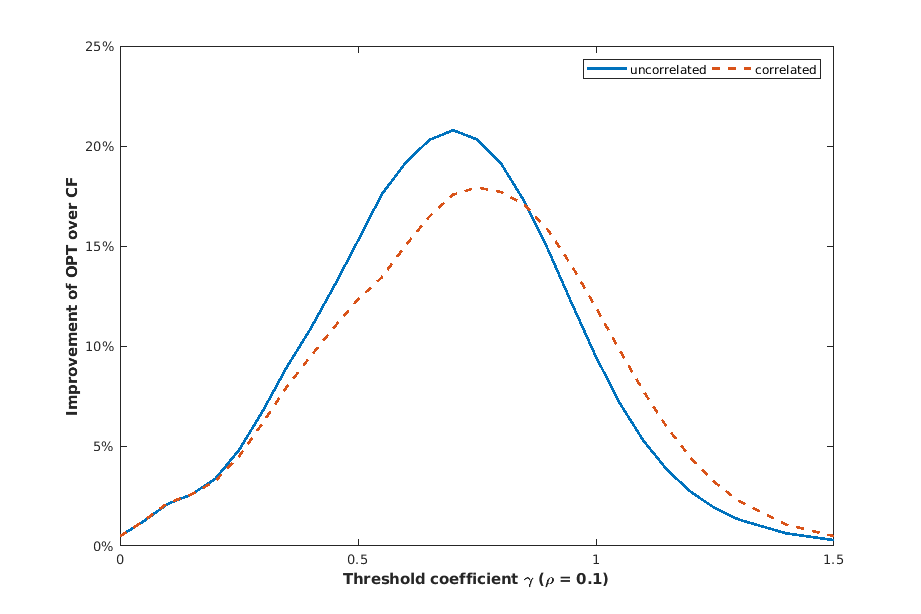}
\caption{graph 4}
\label{fig:varying_thr_4}
\end{subfigure}
\caption{$\delta^{OPT}$ as a function of $\gamma$ for 4 random graphs ($I = 5$, $d = 7$, $\rho = 0.1$)}
\label{fig:varying_thr_coeff}
\end{figure}

The influence of the time threshold $t^*$ (through the parameter $\gamma$) is studied in Figure~\ref{fig:varying_thr_coeff}. Four arbitrary random graphs are chosen, and for each the relative improvement is plotted against $\gamma$, with $\rho = 0.1$. We again observe that significant gains can be made compared to the closest-first policy, and that the scope of this improvement depends on the parameters. Here we can see that the behaviour is similar in both the correlated and uncorrelated cases. If $\gamma$ is close to zero (and hence $t^*$ is too), the OPT policy cannot improve upon the CF policy. The time threshold is too low to meet unless the location of a fire coincides with the location of one of the idle trucks. As a result, the fraction of late arrivals is close to 1 independent of which trucks are sent. As $\gamma$ grows, there is more room for improvement. However, when $\gamma$ approaches 1, the relative improvement of OPT drops to zero again. The reason is that in this case the time threshold $t^*$ is so large it can always be met, even if the dispatching is far from optimal.

\begin{table}[t]
  \caption{Minimum, maximum and mean $\delta^{OPT}$ evaluated over 150 random graphs ($\rho = 0.1$, $\gamma = 0.6$)}
  \label{tbl:aggregate_results}
	\centering
  \includegraphics[width=\linewidth, trim={3.0cm 2.0cm 3.0cm 2.0cm}, clip]{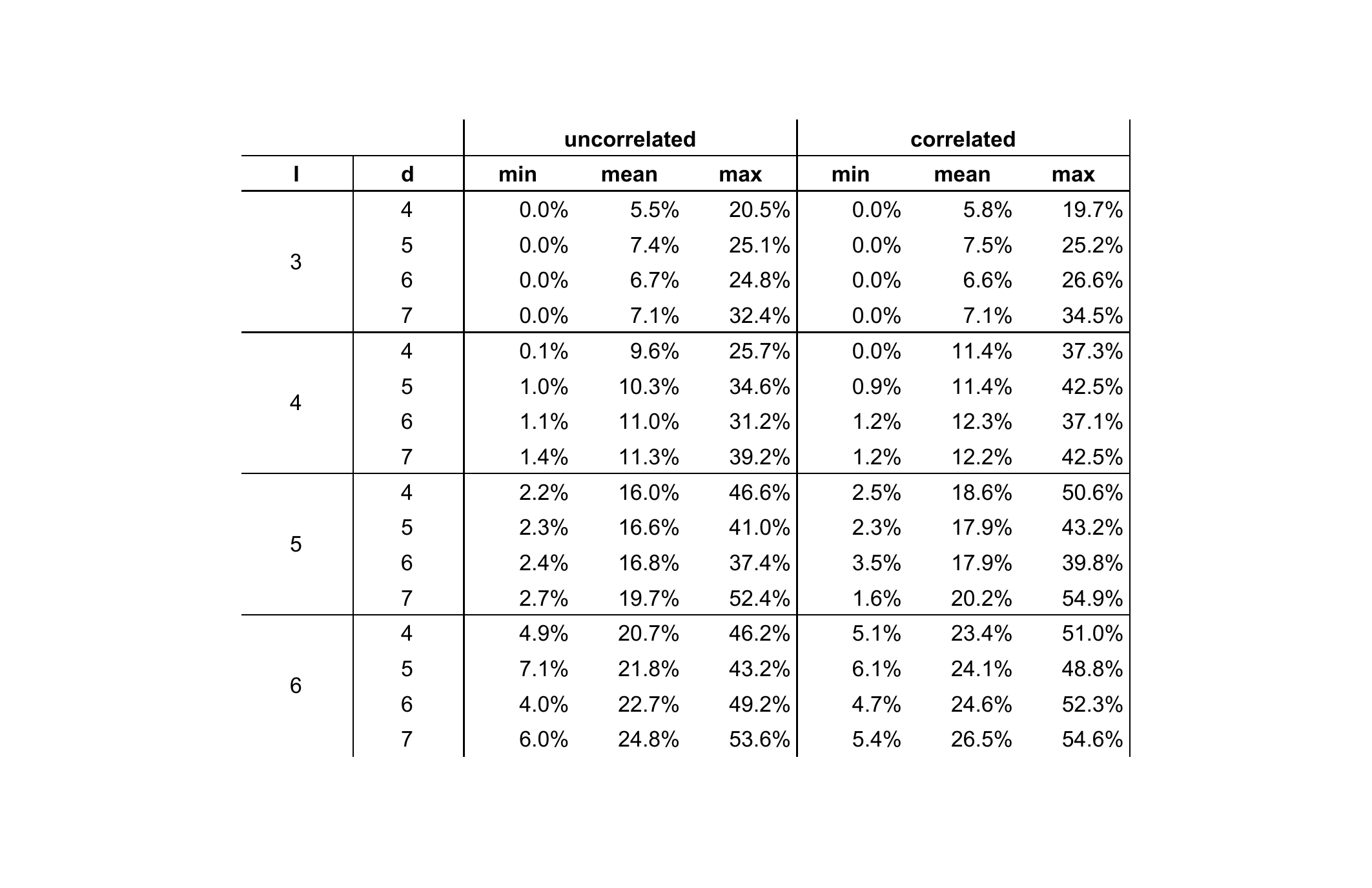}
\end{table}

For a more thorough review of the relative improvement of OPT over CF we turn to Table \ref{tbl:aggregate_results}. This shows the relative improvements a function of the graph size parameter $d$ and the number of trucks $I$, for $\rho=0.1$ and $\gamma=0.6$. For every combination of $I$ and $d$, we generate 150 random graphs. The values in Table~\ref{tbl:aggregate_results} represent the minimum, mean and maximum over these 150 random graphs for each parameter set. We can see a modest increase in relative improvement in $d$, and a significant improvement in $I$, reaching an average improvement of over $20\%$ with $I=6$ trucks, and over $50\%$ for certain instances with driving time correlation.

\begin{figure}[t]
	\centering
	\includegraphics[width=\linewidth]{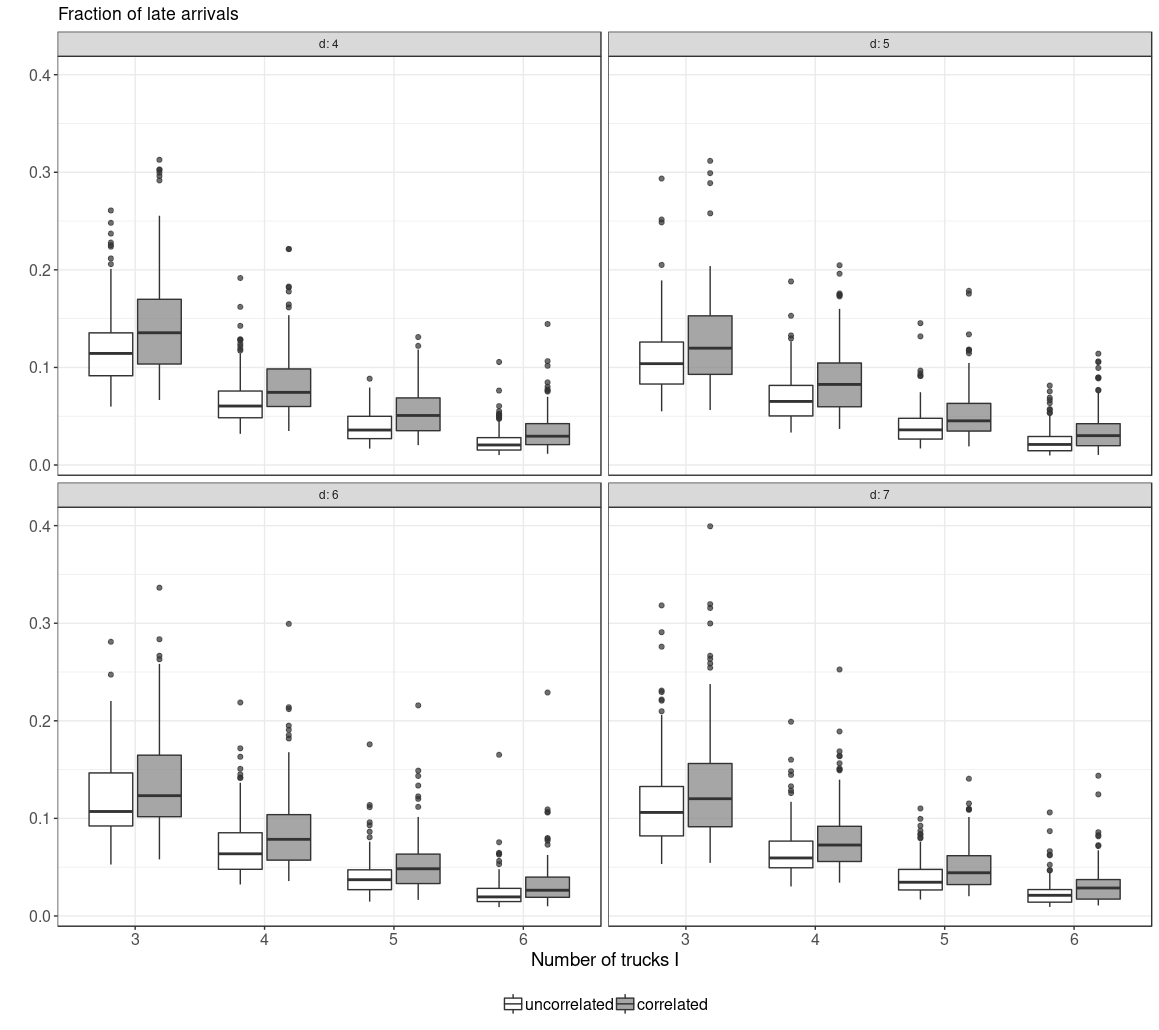}
	\caption{Confidence intervals for performance of the OPT policy for different values of $d$ ($\rho = 0.1$, $\gamma = 0.6$)}
    \label{fig:varying_size_abs}
\end{figure}

In Figure~\ref{fig:varying_size_abs} we show the fraction of late arrivals for OPT for the same set of experiments discussed above. That is, for different values of $d$ and $I$ we plot the confidence interval over all 150 graphs considered. Although we observed from Table \ref{tbl:aggregate_results} that the average relative improvement of OPT over CF is not significantly affected by whether we consider driving-time correlation, we see from Figure \ref{fig:varying_size_abs} that the fraction of late arrivals increases when correlation is taken into account. This indicates that in this case it is more important to deviate from the CF policy in order to limit the fraction of late arrivals. Since in practice there is always some degree of driving-time correlation, these results suggest that when dispatching multiple trucks it is valuable to deviate from CF dispatching. This is in contrast to the case with a single truck, when CF is close to optimal~\cite{jagtenberg2017dynamic}.

\begin{figure}
\begin{subfigure}{.45\textwidth}
\centering
\includegraphics[width=\linewidth, trim={5.5cm 10.0cm 5.0cm 9.0cm}, clip]{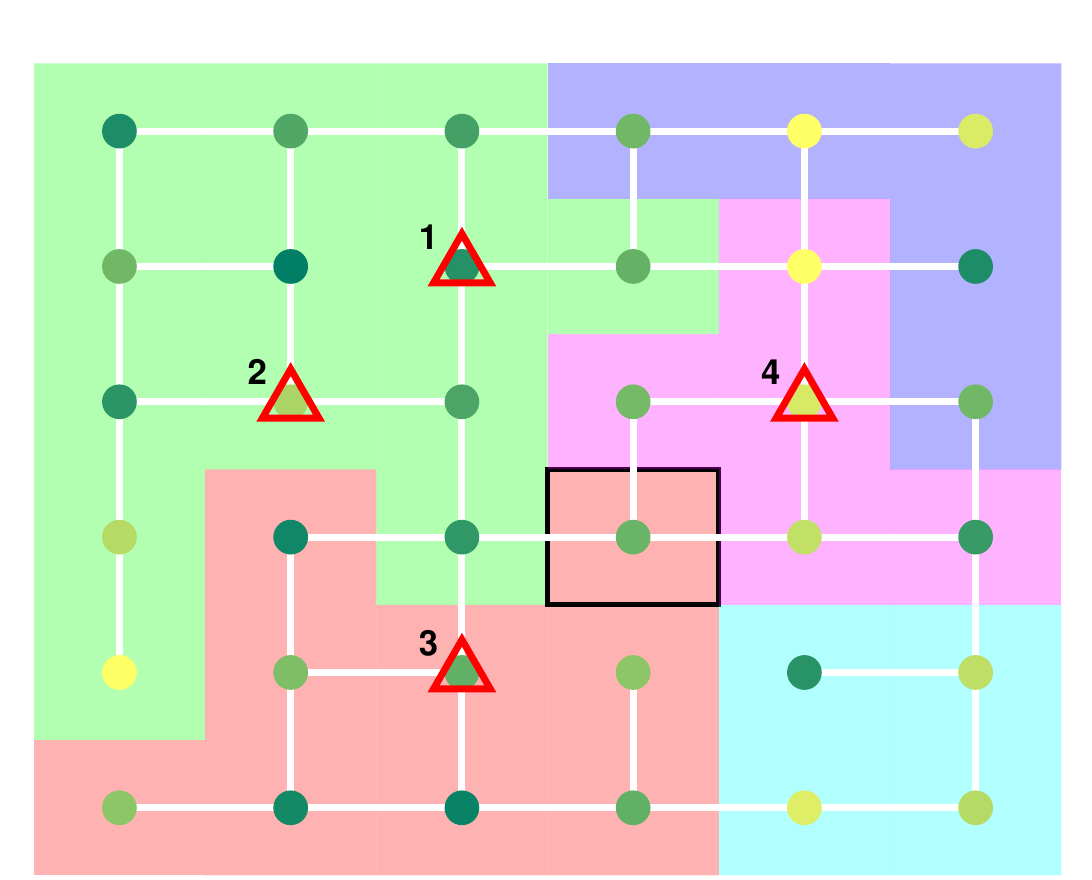}
\caption{$\boldsymbol a^{OPT}_{uc}(\boldsymbol C,j)$}
\label{fig:cor_vs_uncor_0}
\end{subfigure}\hfill
\begin{subfigure}{.05\textwidth}
\centering
\includegraphics[width=\linewidth, height=5cm, trim={10.5cm 9.5cm 9.8cm 9.0cm}, clip]{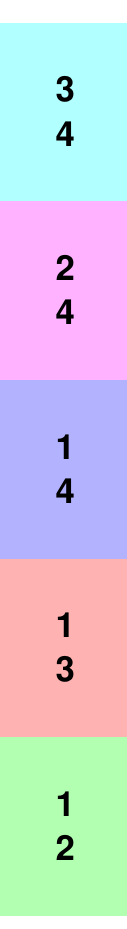}
\end{subfigure}\hfill
\begin{subfigure}{.45\textwidth}
\centering
\includegraphics[width=\linewidth, trim={5.5cm 10.0cm 5.0cm 9.0cm}, clip]{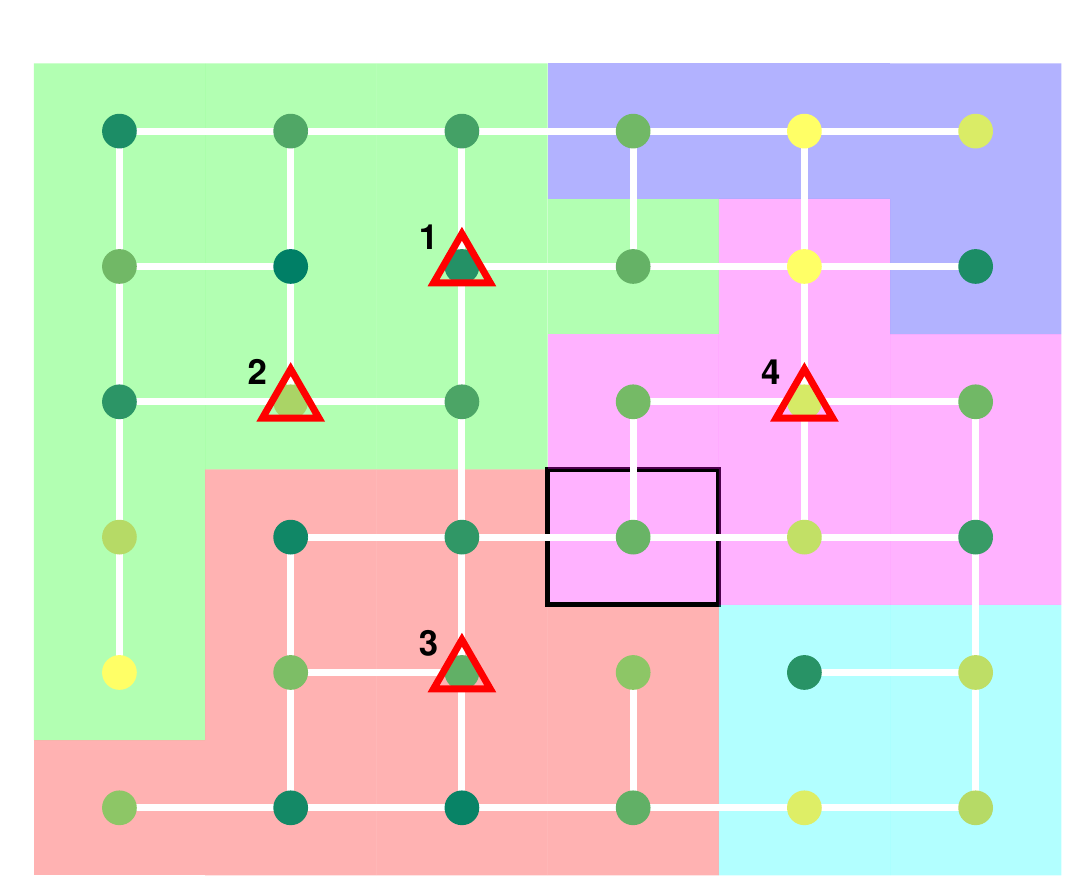}
\caption{$\boldsymbol a^{OPT}_{c}(\boldsymbol C,j)$}
\label{fig:cor_vs_uncor_1}
\end{subfigure}
\caption{Example of difference between $\boldsymbol a^{OPT}_{uc}(\bff,j)$ and $\boldsymbol a^{OPT}_{c}(\bff,j)$ on a random graph}
\label{fig:cor_vs_uncor}
\end{figure}

{\bf Impact of correlation on the optimal policy.} To illustrate the difference between the optimal policies without correlation ($\boldsymbol a^{OPT}_{uc}$) and with correlation ($\boldsymbol a^{OPT}_{c}$) we select a random graph with $d=6$ ($J=36$ demand locations) and $I=4$ trucks, see Figure~\ref{fig:cor_vs_uncor}. The demand locations are coloured according to the arrival rates of new incidents, with green corresponding to low rates. We are looking at the state $\bff = \boldsymbol C$ with all 4 trucks available. The background of each location $j$ is colored according to the corresponding policy $\boldsymbol a^{OPT}(\boldsymbol C,j)$. For example, if a new incident happens at a demand location with green background, then the trucks 1 and 2 will be dispatched. 

While for this particular choice of graph and parameters the impact of correlation is relatively small, it is useful for illustrating how the optimal policy changes when correlation is introduced. For instance, to the demand location highlighted in black in the middle of the graph the policy $\boldsymbol a^{OPT}_{uc}$ dispatches the trucks 1 and 3 that share one edge on their way to that location. The policy $\boldsymbol a^{OPT}_{c}$ instead dispatches the trucks 2 and 4 that share no edges in their shortest paths, as shared edges imply higher probability of being late in the presence of driving-time correlation. 

The other two changes in this example, as well as those in other instances we evaluated, follow a similar pattern: the optimal policy with correlation may be different from the optimal policy without correlation for those demand locations where $\boldsymbol a^{OPT}_{uc}$ dispatches two trucks with overlapping routes. However, this need not be the case, and the example in Figure~\ref{fig:cor_vs_uncor} also includes such demand locations where $\boldsymbol a^{OPT}_{c}$ remains unchanged compared to $\boldsymbol a^{OPT}_{uc}$, because in these cases the decrease in expected response time when changing actions does not outweigh the coverage reduction resulting from this change. This illustrates the complexity of finding the optimal policy for this model, and the difficulties one would encounter when trying to generalize the observations obtained from Figure~\ref{fig:cor_vs_uncor} into some kind of heuristic. One main reason for this is the complex interactions encountered in this model. For instance, changing the arrival rate in one part of the network may affect the optimal policy elsewhere.

To see the extend to which driving-time correlation affects the optimal policy for a broader range of instances we conduct the following experiment. We generate 150 random graphs, and for every graph we compute $\boldsymbol a^{OPT}_{c}$ and $\boldsymbol a^{OPT}_{uc}$. In order to study the impact of ignoring driving-time correlation, we look what happens with the system performance if we use $\boldsymbol a^{OPT}_{uc}$ in a setting with driving-time correlation. To do this we plug the policy $\boldsymbol a^{OPT}_{uc}$ into the Bellman equations~\eqref{eq:bellman_sd} for a fixed policy with the costs corresponding to the correlated case, and measure the relative increase in value function compared to the policy $\boldsymbol a^{OPT}_{c}(\bff,j)$. Note that the relative increase in value function is equivalent to the relative increase in the fraction of late arrivals.

\begin{table}[]
\caption{Relative increase in fraction of late arrivals when neglecting correlation ($d = 6$, $\rho = 0.1$, $\gamma = 0.6$)}
\label{tbl:cor_vs_uncor}
\centering
\begin{tabular}{l|l|l|l}
\multicolumn{1}{c|}{\textbf{I}} & \multicolumn{1}{c|}{\textbf{min}} & \multicolumn{1}{c|}{\textbf{mean}} & \multicolumn{1}{c}{\textbf{max}} \\ \hline
3                               & 0.0\%                             & 1.3\%                              & 7.4\%                            \\ \hline
4                               & 0.2\%                             & 2.8\%                              & 12.2\%                           \\ \hline
5                               & 0.3\%                             & 4.8\%                              & 16.3\%                           \\ \hline
6                               & 1.1\%                             & 7.1\%                              & 21.3\%                          
\end{tabular}
\end{table}

Table~\ref{tbl:cor_vs_uncor} shows the aggregate results of this experiment with minimum, maximum and mean relative increase in fraction of late arrivals computed over 150 random graphs. We see that the importance of taking driving-time correlation into account grows with the number of trucks in the system. With more vehicles available there are more options for making a dispatching decision to avoid potential traffic jams for the current and upcoming incidents. The average decrease in performance when using the policy derived under the assumption of uncorrelated driving times in a setting with driving-time correlation reached $7.1\%$ for 6 trucks, and for some instances was over $20\%$.

\subsection{Performance of the heuristics}\label{sec:heur}

{\bf Improvement over closest-first.}  In this section we compare the performance of the two heuristics OSI and OSIA to the optimal policy OPT, both in terms of fraction of late arrivals and computational time.

\renewcommand{\arraystretch}{0.8}
\begin{table}[h]
\footnotesize
\caption{Aggregate performance evaluated over multiple random graphs ($d=6$, $\gamma = 0.6$)}
\label{tbl:aggregate_results_expand}
\centering
\begin{tabular}{ll|l|l|l|l|l|l|l|l|}
\multicolumn{1}{c}{\textbf{}}           & \multicolumn{1}{c|}{\textbf{}}                    & \multicolumn{4}{c|}{\textbf{uncorrelated}}                                                                                                                                                                                                                                                               & \multicolumn{4}{c|}{\textbf{correlated}}                                                                                                                                                                                                                                                                 \\ \hline
\multicolumn{1}{c|}{\pmb{$I$}}         & \multicolumn{1}{c|}{\pmb{$\rho$}} & \multicolumn{1}{c|}{\pmb{$FLAR^{CF}$}} & \multicolumn{1}{c|}{\pmb{$\delta^{OPT}$}} & \multicolumn{1}{c|}{\pmb{$\delta^{OSI}$}} & \multicolumn{1}{c|}{\pmb{$\delta^{OSIA}$}} & \multicolumn{1}{c|}{\pmb{$FLAR^{CF}$}} & \multicolumn{1}{c|}{\pmb{$\delta^{OPT}$}} & \multicolumn{1}{c|}{\pmb{$\delta^{OSI}$}} & \multicolumn{1}{c|}{\pmb{$\delta^{OSIA}$}} \\ \hline
\multicolumn{1}{l|}{\multirow{5}{*}{3}} & 0.02                                              & 0.39\%                                                     & 4.83\%                                                                       & 4.83\%                                                                       & 2.81\%                                                                        & 0.51\%                                                     & 5.54\%                                                                       & 5.54\%                                                                       & 4.12\%                                                                        \\
\multicolumn{1}{l|}{}                   & 0.04                                              & 0.48\%                                                     & 6.36\%                                                                       & 6.36\%                                                                       & 4.88\%                                                                        & 0.60\%                                                     & 6.52\%                                                                       & 6.52\%                                                                       & 5.43\%                                                                        \\
\multicolumn{1}{l|}{}                   & 0.1                                               & 0.77\%                                                     & 6.70\%                                                                       & 6.70\%                                                                       & 6.36\%                                                                        & 0.89\%                                                     & 6.60\%                                                                       & 6.60\%                                                                       & 6.33\%                                                                        \\
\multicolumn{1}{l|}{}                   & 0.4                                               & 2.12\%                                                     & 2.57\%                                                                       & 2.57\%                                                                       & 2.51\%                                                                        & 2.19\%                                                     & 2.59\%                                                                       & 2.59\%                                                                       & 2.52\%                                                                        \\
\multicolumn{1}{l|}{}                   & 0.6                                               & 2.74\%                                                     & 1.48\%                                                                       & 1.48\%                                                                       & 1.45\%                                                                        & 2.79\%                                                     & 1.50\%                                                                       & 1.50\%                                                                       & 1.47\%                                                                        \\ \hline
\multicolumn{1}{l|}{\multirow{5}{*}{4}} & 0.02                                              & 0.20\%                                                     & 9.49\%                                                                       & 9.49\%                                                                       & 5.10\%                                                                        & 0.31\%                                                     & 12.24\%                                                                      & 12.24\%                                                                      & 9.53\%                                                                        \\
\multicolumn{1}{l|}{}                   & 0.04                                              & 0.25\%                                                     & 11.55\%                                                                      & 11.54\%                                                                      & 8.51\%                                                                        & 0.37\%                                                     & 13.53\%                                                                      & 13.52\%                                                                      & 11.43\%                                                                       \\
\multicolumn{1}{l|}{}                   & 0.1                                               & 0.47\%                                                     & 11.02\%                                                                      & 10.99\%                                                                      & 10.45\%                                                                       & 0.59\%                                                     & 12.32\%                                                                      & 12.29\%                                                                      & 11.87\%                                                                       \\
\multicolumn{1}{l|}{}                   & 0.4                                               & 1.82\%                                                     & 3.77\%                                                                       & 3.75\%                                                                       & 3.28\%                                                                        & 1.90\%                                                     & 4.13\%                                                                       & 4.12\%                                                                       & 3.50\%                                                                        \\
\multicolumn{1}{l|}{}                   & 0.6                                               & 2.50\%                                                     & 2.14\%                                                                       & 2.14\%                                                                       & 1.94\%                                                                        & 2.56\%                                                     & 2.33\%                                                                       & 2.32\%                                                                       & 2.05\%                                                                        \\ \hline
\multicolumn{1}{l|}{\multirow{5}{*}{5}} & 0.02                                              & 0.10\%                                                     & 15.71\%                                                                      & 15.62\%                                                                      & 8.56\%                                                                        & 0.18\%                                                     & 17.38\%                                                                      & 17.30\%                                                                      & 13.50\%                                                                       \\
\multicolumn{1}{l|}{}                   & 0.04                                              & 0.14\%                                                     & 18.57\%                                                                      & 18.35\%                                                                      & 13.60\%                                                                       & 0.21\%                                                     & 19.56\%                                                                      & 19.32\%                                                                      & 16.22\%                                                                       \\
\multicolumn{1}{l|}{}                   & 0.1                                               & 0.29\%                                                     & 16.84\%                                                                      & 16.49\%                                                                      & 14.35\%                                                                       & 0.38\%                                                     & 17.94\%                                                                      & 17.58\%                                                                      & 16.00\%                                                                       \\
\multicolumn{1}{l|}{}                   & 0.4                                               & 1.60\%                                                     & 4.76\%                                                                       & 4.68\%                                                                       & 3.90\%                                                                        & 1.68\%                                                     & 5.22\%                                                                       & 5.12\%                                                                       & 4.06\%                                                                        \\
\multicolumn{1}{l|}{}                   & 0.6                                               & 2.34\%                                                     & 2.56\%                                                                       & 2.53\%                                                                       & 2.21\%                                                                        & 2.40\%                                                     & 2.79\%                                                                       & 2.75\%                                                                       & 2.32\%                                                                        \\ \hline
\multicolumn{1}{l|}{\multirow{5}{*}{6}} & 0.02                                              & 0.05\%                                                     & 20.45\%                                                                      & 20.15\%                                                                      & 11.73\%                                                                       & 0.11\%                                                     & 22.00\%                                                                      & 21.70\%                                                                      & 17.33\%                                                                       \\
\multicolumn{1}{l|}{}                   & 0.04                                              & 0.07\%                                                     & 24.90\%                                                                      & 24.37\%                                                                      & 17.89\%                                                                       & 0.13\%                                                     & 25.91\%                                                                      & 25.30\%                                                                      & 21.64\%                                                                       \\
\multicolumn{1}{l|}{}                   & 0.1                                               & 0.18\%                                                     & 22.75\%                                                                      & 22.05\%                                                                      & 18.44\%                                                                       & 0.25\%                                                     & 24.63\%                                                                      & 23.75\%                                                                      & 21.57\%                                                                       \\
\multicolumn{1}{l|}{}                   & 0.4                                               & 1.43\%                                                     & 5.99\%                                                                       & 5.84\%                                                                       & 4.92\%                                                                        & 1.50\%                                                     & 6.67\%                                                                       & 6.49\%                                                                       & 5.17\%                                                                        \\
\multicolumn{1}{l|}{}                   & 0.6                                               & 2.20\%                                                     & 3.16\%                                                                       & 3.10\%                                                                       & 2.70\%                                                                        & 2.26\%                                                     & 3.47\%                                                                       & 3.40\%                                                                       & 2.87\%                                                                       
\end{tabular}
\end{table}

Table \ref{tbl:aggregate_results_expand} shows the relative difference of OSI and OSIA with CF, in addition to that of OPT. The values of $\delta^{OSI}$ and $\delta^{OSIA}$ are computed the same way as $\delta^{OPT}$. The numbers presented in the table are the mean values of the corresponding metrics evaluated over 150 randomly generated graphs. The values of $d$ and $\gamma$ are fixed, and we vary the load $\rho$ (in the range $\{0.02, 0.04, 0.1, 0.4, 0.6\}$) and the number of trucks $I$ (in the range 3 to 6). For every combination of $I$ and $\rho$ the minimum, mean and maximum over 150 randomly generated graphs is presented. The improvement over CF for all three policies first increases with $\rho$ followed by a decrease for high loads. Both OSI and OSIA policies show significant improvement over the CF policy for lower values of $\rho$, and are relatively close to the performance of OPT. The improvement over CF grows with $I$ and is larger in the presence of driving-time correlation, similar to what we observed in Table~\ref{tbl:aggregate_results}. As it can be seen from the more detailed Tables~\ref{tbl:aggregate_results_expand2}, \ref{tbl:aggregate_results_expand3} and \ref{tbl:aggregate_results_expand5} in Appendix~\ref{sec:add_num}, the heuristics performance also improves as $d$ increases, suggesting that their performance is better for larger networks. Appendix~\ref{sec:add_num} also includes Table~\ref{tbl:aggregate_results_expand4}, which shows the relative improvement of OSIA over CF for $\rho = 0.02$ and $I=7$ for larger values of $d$. Here we see that as $I$ and $d$ grow larger, the gap with CF increases as well.

\rev{Note that, in our setting, the fraction of late arrivals under the CF policy $FLAR^{CF}$ is relatively low. This would mean that the improvement over CF in the number of late arrivals is low compared to the total number of incidents. However, this improvement should not be understated. From the emergency services perspective, any improvement in late arrivals is considered significant. In the case of FDAA, for example, the original idea of dispatching two trucks instead of one is targeted at reducing the risk of a possible delay, despite additional operational costs. This shows the importance of any decrease in response time. The further gains that can be achieved by changing the dispatching strategy are particularly valuable, given that it does not involve any extra operational costs.}

\renewcommand{\arraystretch}{0.8}
\begin{table}[h]
\footnotesize
\caption{Average optimality gap of OSI and OSIA ($d = 6$, $\gamma = 0.6$)}
\label{tbl:optimality_gap}
\centering
\begin{tabular}{ll|l|l|l|l|}
\multicolumn{1}{c}{\textbf{}}           & \multicolumn{1}{c|}{\textbf{}}                    & \multicolumn{2}{c|}{\textbf{uncorrelated}}                                                                                                                   & \multicolumn{2}{c|}{\textbf{correlated}}                                                                                                                     \\ \hline
\multicolumn{1}{c|}{\pmb{I}}         & \multicolumn{1}{c|}{\pmb{$\rho$}} & \multicolumn{1}{c|}{\textbf{OSI}} & \multicolumn{1}{c|}{\textbf{OSIA}} & \multicolumn{1}{c|}{\textbf{OSI}} & \multicolumn{1}{c|}{\textbf{OSIA}} \\ \hline
\multicolumn{1}{l|}{\multirow{5}{*}{3}} & 0.02                                              & 0.00\%                                                                       & 2.36\%                                                                        & 0.00\%                                                                       & 1.68\%                                                                        \\
\multicolumn{1}{l|}{}                   & 0.04                                              & 0.00\%                                                                       & 1.72\%                                                                        & 0.00\%                                                                       & 1.24\%                                                                        \\
\multicolumn{1}{l|}{}                   & 0.1                                               & 0.00\%                                                                       & 0.37\%                                                                        & 0.00\%                                                                       & 0.29\%                                                                        \\
\multicolumn{1}{l|}{}                   & 0.4                                               & 0.00\%                                                                       & 0.06\%                                                                        & 0.00\%                                                                       & 0.07\%                                                                        \\
\multicolumn{1}{l|}{}                   & 0.6                                               & 0.00\%                                                                       & 0.03\%                                                                        & 0.00\%                                                                       & 0.03\%                                                                        \\ \hline
\multicolumn{1}{l|}{\multirow{5}{*}{4}} & 0.02                                              & 0.00\%                                                                       & 6.01\%                                                                        & 0.01\%                                                                       & 3.57\%                                                                        \\
\multicolumn{1}{l|}{}                   & 0.04                                              & 0.01\%                                                                       & 3.91\%                                                                        & 0.01\%                                                                       & 2.70\%                                                                        \\
\multicolumn{1}{l|}{}                   & 0.1                                               & 0.03\%                                                                       & 0.66\%                                                                        & 0.04\%                                                                       & 0.52\%                                                                        \\
\multicolumn{1}{l|}{}                   & 0.4                                               & 0.01\%                                                                       & 0.51\%                                                                        & 0.02\%                                                                       & 0.66\%                                                                        \\
\multicolumn{1}{l|}{}                   & 0.6                                               & 0.01\%                                                                       & 0.21\%                                                                        & 0.01\%                                                                       & 0.28\%                                                                        \\ \hline
\multicolumn{1}{l|}{\multirow{5}{*}{5}} & 0.02                                              & 0.18\%                                                                       & 11.67\%                                                                       & 0.16\%                                                                       & 5.98\%                                                                        \\
\multicolumn{1}{l|}{}                   & 0.04                                              & 0.36\%                                                                       & 7.24\%                                                                        & 0.40\%                                                                       & 4.79\%                                                                        \\
\multicolumn{1}{l|}{}                   & 0.1                                               & 0.47\%                                                                       & 3.05\%                                                                        & 0.49\%                                                                       & 2.41\%                                                                        \\
\multicolumn{1}{l|}{}                   & 0.4                                               & 0.09\%                                                                       & 0.90\%                                                                        & 0.11\%                                                                       & 1.22\%                                                                        \\
\multicolumn{1}{l|}{}                   & 0.6                                               & 0.04\%                                                                       & 0.36\%                                                                        & 0.04\%                                                                       & 0.49\%                                                                        \\ \hline
\multicolumn{1}{l|}{\multirow{5}{*}{6}} & 0.02                                              & 0.67\%                                                                       & 16.26\%                                                                       & 0.58\%                                                                       & 8.13\%                                                                        \\
\multicolumn{1}{l|}{}                   & 0.04                                              & 1.04\%                                                                       & 11.72\%                                                                       & 1.12\%                                                                       & 7.04\%                                                                        \\
\multicolumn{1}{l|}{}                   & 0.1                                               & 1.01\%                                                                       & 5.73\%                                                                        & 1.33\%                                                                       & 4.18\%                                                                        \\
\multicolumn{1}{l|}{}                   & 0.4                                               & 0.16\%                                                                       & 1.13\%                                                                        & 0.20\%                                                                       & 1.60\%                                                                        \\
\multicolumn{1}{l|}{}                   & 0.6                                               & 0.06\%                                                                       & 0.47\%                                                                        & 0.07\%                                                                       & 0.63\%                                                                       
\end{tabular}
\end{table}

Given the value functions $g^{OPT}$ and $g^{OSIA}$ of the OPT and OSIA policies, respectively, we compute the OSIA optimality gap as $\frac{g^{OSIA}-g^{OPT}}{g^{OPT}}\times 100\%$. We compute the optimality gap for OSI in a similar way. 
Table~\ref{tbl:optimality_gap} shows the average optimality gap of the OSI and OSIA policies computed over 150 random graphs for each combination of $I$ and $\rho$. The performance of both OSI and OSIA stays within a few percent of OPT. The optimality gap grows with $I$. The OSI policy performs slightly better in a setting without correlation, while the opposite is true for OSIA. The optimality gap of both OSI and OSIA decreases in $\rho$, suggesting that these approximations perform best in the high load regime. Note that while the optimality gap of these heuristics grows in the network size, we have seen from Tables~\ref{tbl:aggregate_results_expand}, \ref{tbl:aggregate_results_expand2}, \ref{tbl:aggregate_results_expand3}, \ref{tbl:aggregate_results_expand5} and~\ref{tbl:aggregate_results_expand4} that the improvement over CF also does. So while neither OSI nor OSIA is asymptotically optimal, their performance in fact improves as the network grows larger.

{\bf Computational time.} Next, we take a look at the computational time of the various policies. If  computational time would not be an issue, then using the OPT policy is an obvious choice. However, solving MDP exactly quickly becomes problematic when the instance size grows, as the size of the state space grows exponentially in $I$. In our experiments, the main issue with solving the MDP exactly for larger instances was not the running time of policy iteration, but the size of the array with transition probabilities ($|\mathcal{S}|\times|\mathcal{S}|\times|\mathcal{A}|$). As a result, computing the OPT policy breaks down for even moderate-sized networks (e.g., $I=7$, $d=6$).

\begin{figure}
	\centering
	\includegraphics[width=14cm]{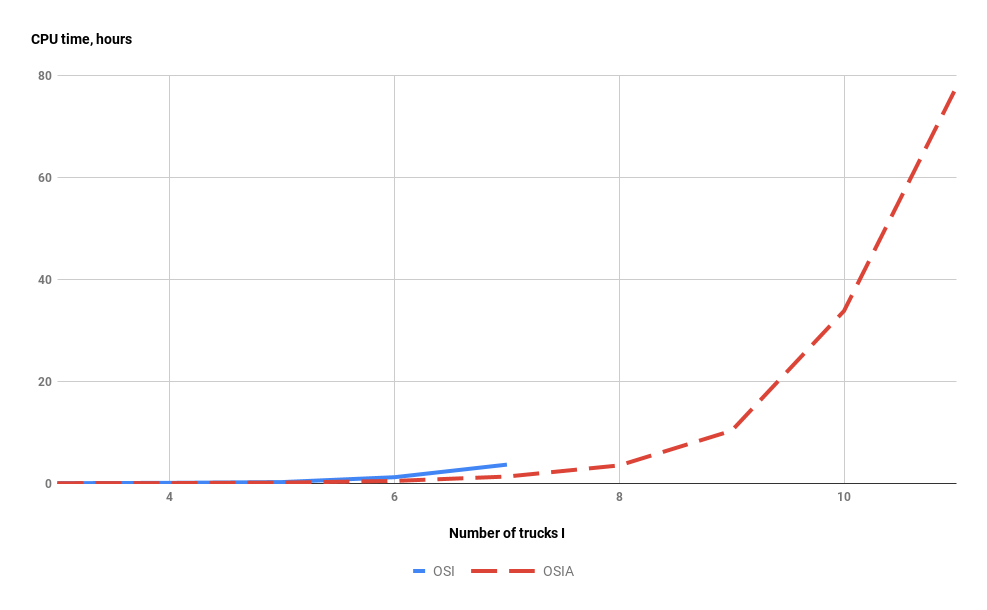}
	\caption{Change in computational time as $I$ grows ($J = 625$)}
    \label{fig:varying_I_cpu}
\end{figure}

\begin{figure}
	\centering
	\includegraphics[width=14cm]{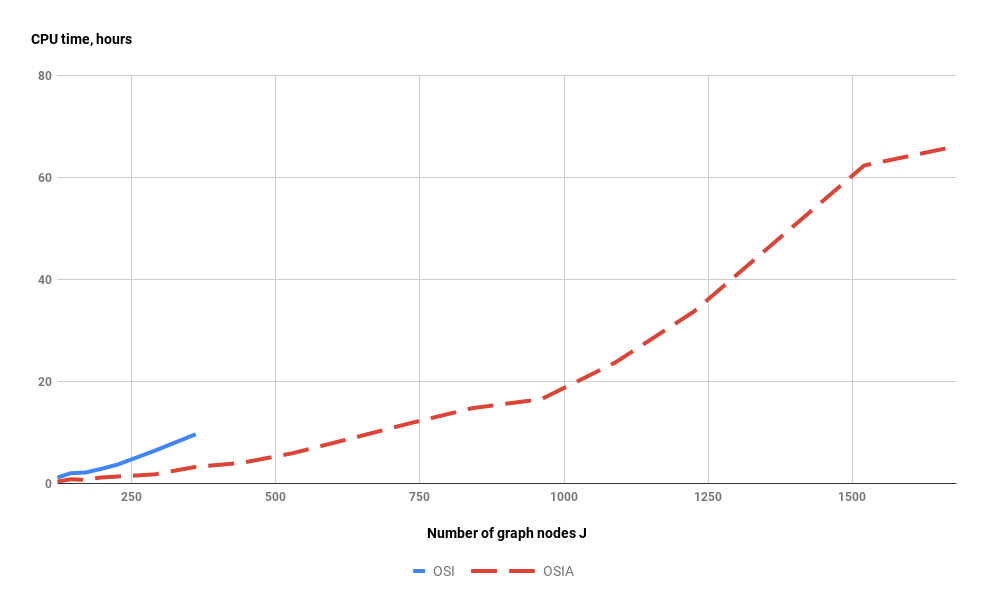}
	\caption{Change in computational time as $J$ grows ($I = 7$)}
    \label{fig:varying_J_cpu}
\end{figure}

To compare the computational performance of OSI and OSIA, we plot the computational time for determining these policies against $I$ (Figure~\ref{fig:varying_I_cpu}) and the number of demand locations $J=d^2$ (Figure~\ref{fig:varying_J_cpu}). Here we use a single randomly generated graph for each data point. The OSI policy is computed faster then the optimal, but still requires solving a set of $|\mathcal{S}+1|$ Bellman equations. Storing a $|\mathcal{S}+1|\times|\mathcal{S}+1|$ matrix of coefficients for the system of Bellman equations becomes infeasible, which is why we can only determine the OSI policy for small values of $I$ and $J$. The computational time of the OSIA policy  shows significantly slower growth in $I$ and $J$ than that of OSI. Moreover, it does not require storing large data structures, and makes it feasible to obtain a good policy for problem instances of realistic size.

The computational time of the OSIA heuristic is reasonable for the systems used in our numerical experiments. The algorithm is meant to be used in the offline regime, only once for a given system, and produces look-up tables indicating the dispatching decision to be made for each state of the system. Moreover, in our experiments we ran approximation Algorithm~\ref{alg:approx} sequentially for each state. In real-life applications the OSIA computational time can be significantly decreased by means of parallelization.

\section{Conclusion}\label{sec:conclusion}

In the present work we studied a dispatching problem in a fire department where two trucks have to be dispatched to an incident location, and the decision is to be made on which idle trucks to send. We modelled the region served by a fire department as a connected graph and formulated the dispatching problem as an MDP. The optimal policy was obtained by solving the MDP exactly using policy iteration. 

Using small problem instances, we showed that the current practice of sending the two closest trucks can be far from optimal, with optimality gap reaching $50\%$ in certain cases. As obtaining the optimal policy for large problem instances is computationally infeasible, we also derived a one-step-improvement OSI policy, that can be obtained faster and for larger problem instances than OPT. In our experiments, however, OSI still remained computationally infeasible for problem instances of realistic sizes. Therefore, we introduced the OSIA policy that incorporates an approximation scheme into the OSI policy computation procedure. The OSIA policy performed close to the optimal performance with optimality gap of about $2\%$, and significantly lower computational time that allows for solving real-life sized problem instances.

We considered two types of stochastic behaviour in driving time when two trucks are dispatched to the same incident location. If two trucks traverse the same edge in a graph we assume their travelling times to be either independent of each other (uncorrelated), or the same (correlated). Our experiments show that introducing correlation makes a difference compared to sending two closest trucks, even if the load is small. Since performance is measured based on response time, sending two closest is not necessarily optimal anymore.


As discussed in Section~\ref{sec:num_opt}, analytically characterizing the optimal policy for general networks seems untractable, due to the complex network dynamics that may propagate even small perturbations throughout the network. However, we are optimistic that for small network instances or specific network structures (such as linear networks), one may be able to obtain structural results on the optimal policy. Doing this for both the case with and without correlation may lead to interesting insights into where and how these two optimal policies differ.

This work can be extended in several interesting ways. For instance, the model and results could be modified to accommodate the following:
\begin{itemize}
\item[-] Instead of only considering perfect or no correlation between the driving time, we could allow for intermediate levels of correlation by assuming that the driving time on a single edge is hyperexponential instead of exponential. By coupling only one of the branches of this distribution we can accommodate any correlation coefficient.
\item[-] Note that changing the driving time distribution does not affect the MDP formulation, but rather the immediate costs. So in order to allow for driving time distributions beyond exponential we would have to generalize Proposition~\ref{pro:Erlang}. Note that if we use a heavy-tailed distribution, the results can potentially show a more significant advantage of using the OPT policy instead of CF. We expect a larger optimality gap for CF in the case of heavy-tailed driving time distribution since the larger variance in response time necessitates more careful dispatching.
\item[-] The MDP formulation itself can be enhanced by allowing more than two trucks to be dispatched to an incident, and we can generalize the definition of the response time accordingly. This would entail changing the action space from all actions that dispatch at most 2 trucks to those that dispatch at most $k$ trucks. The main difficulty in making this extension lies in computing the immediate cost $\mathds{P}(R(\bfa(\boldsymbol f,j),j)>t^*)$ for those actions $\bfa$ that dispatch more than two trucks. If only the first truck to arrive is relevant, the costs can be computed along the lines of Proposition~\ref{pro:Erlang}, by conditioning on the realizations of the driving times of all trucks. If the performance metric depends on more than just the first truck to arrive, generalizing the results obtained here may be more complex.
\item[-] When two trucks are dispatched from the same station, we may assume that each takes a different path in order to avoid driving-time correlation. Including this in the model may result in the optimal policy and heuristics to dispatch trucks from the same station more often.
\end{itemize}

\noindent
\textbf{Acknowledgements}

This research was funded by an NWO grant, under contract number 438-15-506. We would like to thank Guido Legemaate from the Fire Department Amsterdam-Amstelland for the useful discussions and insights that provided motivation for this paper. We also like to thank the anonymous referees, whose suggestions have led to a significant improvement of the paper.

\bibliographystyle{plain}
\bibliography{references}

\newpage
\appendix

\nomenclature[01]{$\mathcal{J}=\{1,...,J\}$}{Set of demand location / nodes in a graph}
\nomenclature[02]{$E$}{Set of edges in a graph}
\nomenclature[03]{$\mathcal{I} \subseteq \mathcal{J}$}{Demand locations containing a fire station}
\nomenclature[04]{$I = \lvert \mathcal{I} \rvert$}{Number of fire stations}
\nomenclature[05]{$C_i$}{Number of fire trucks with the base station $i \in \mathcal{I}$}
\nomenclature[06]{$\lambda_j$}{Arrival rate of new fires at a location $j \in \mathcal{J}$}
\nomenclature[07]{$1/\mu$}{Expected time a truck remains busy after being dispatched}
\nomenclature[08]{$\rho = \frac{\sum_{j \in \mathcal{J}}\lambda_j}{I\mu}$}{Load of the system, that is, the amount of work per fire truck per time unit}
\nomenclature[09]{$f_i$}{Number of idle trucks at a station $i \in \mathcal{I}$}
\nomenclature[10]{$\boldsymbol e_i$}{vector of length $I$ with $i$th element equal to 1, and all other elements equal to zero}
\nomenclature[11]{$\bff = (f_1,...,f_I)$}{Vector representing the state of the system}
\nomenclature[12]{$\bfa(\bff, j) = (a_1(\bff, j),...,a_I(\bff, j))$}{The dispatch action taken if a new fire starts at a location $j$ when in state $\bff$}
\nomenclature[13]{$0 \leq a_i(\bff, j) \leq f_i$}{Number of trucks dispatched from station $i \in \mathcal{I}$}
\nomenclature[14]{$\mathcal{S} = \{(f_1,...,f_I) \vert 0 \leq f_i \leq C_i \ \forall i \in \mathcal{I}\}$}{The system state space}
\nomenclature[15]{$s(i, j)$}{Shortest path between nodes $i$ and $j$ in a graph}
\nomenclature[16]{$T_{i,j} = \sum_{e \in s(i,j)} X_e$}{Traveling time between nodes $i$ and $j$, where $X_e \sim {\rm exp}(1)$}
\nomenclature[17]{$T_0$}{Traveling time from a neighboring region to any demand location}
\nomenclature[18]{$R(\bfa,j)$}{Response time to a fire at a location $j$ given a dispatch decision $\bfa$}
\nomenclature[19]{$\tau = \sum_{j\in \mathcal{J}}\lambda_j + \mu \sum_{i \in \mathcal{I}}C_i$}{Transition rate out of any state}
\nomenclature[20]{$\mathcal{A}(\bff)$}{Actions space in state $\bff \in \mathcal{S}$}
\nomenclature[21]{$g^*$}{Average cost incurred per time unit}
\nomenclature[22]{$h^*(\bff)$}{Relative cost incurred over infinite time horizon when starting in state $\bff \in \mathcal{S}$ compared to paying $g^*$ every time unit}
\nomenclature[23]{$\sigma_j(k) \in \mathcal{I}$}{Fire station that is the base station for the $k$th closest truck to location $j$, assuming that truck is idle}
\nomenclature[24]{$k_i, i = 1, 2$}{Number of the closest and the second-closest idle truck in the list $\sigma_j(k)$, $j\in\mathcal{J}, k \in \{1,\dots,\sum_i C_i\}$}
\nomenclature[25]{$J(\bff, t)$}{Expected total cost under the CF policy during the time interval $[0,t]$ starting from state $\bff$}
\nomenclature[26]{$T$}{Parameter of the OSIA heuristic indicating the time it takes for the system to get into the steady state by assumption}
\nomenclature[27]{$D_i$}{Arrival rate of requests for the truck at station $i$}
\nomenclature[28]{$\rho_i  = D_i/\mu$}{Load of the $M/M/1/1$ queue representing fire station $i$}
\nomenclature[29]{$p_i$}{Busy probability of station $i$}
\nomenclature[30]{$t^*$}{Response time threshold}
\nomenclature[31]{$\gamma$}{Parameter that defines the response time threshold $t^*$ for a given graph as a fraction of the maximum traveling time between two nodes}

\renewcommand{\nomname}{List of Notations}
\printnomenclature[2.5in]

\section{Proof of Proposition~\ref{pro:Erlang}}\label{app:proof}

\begin{proof}
The first statement can be readily proven by using the independence of $Y_1$ and $Y_2$:
$$
P(\min\{Y_1, Y_2\} > t^*) = P(Y_1 \geq t^*)P(Y_2 \geq t^*).
$$
Substituting in the distribution of $Y_1$ and $Y_2$ we obtain the desired result.

For the second statement we condition on the value of $Y_0$ to obtain the following expression:
\begin{align*}
  P(Y_0 + \min\{Y_1, Y_2\} > t^*)
  &= \int_{y_0 = 0}^{\infty}f_{Y_0}(y_0)P(\min\{Y_1, Y_2\} > t^* - y_0){\rm d}y_0  \\
  &= \int_{y_0 = 0}^{t^*}f_{Y_0}(y_0)P(Y_1 > t^* - y_0)P(Y_2 > t^* - y_0){\rm d}y_0 + \int_{y_0 = t^*}^{\infty}f_{Y_0}(y_0){\rm d}y_0.
  \end{align*}
 By substituting the distribution function of $Y_0$, $Y_1$ and $Y_2$, and exchanging the order of integration and summation we obtain
 \begin{align*}
  P(R(\bfa,j)>t^*) &= \int_{y_0 = 0}^{t^*}f_{Y_0}(y_0)\sum_{n=0}^{w_1-1}\frac{(t^*-y_0)^n}{n!}e^{-t^*+y_0}\sum_{m=0}^{w_2-1}\frac{(t^*-y_0)^m}{m!}e^{-t^*+y_0}{\rm d}y_0 + \sum_{n=0}^{w_0-1}\frac{{t^*}^n}{n!}e^{-t^*}  \\
  &=\sum_{n=0}^{w_1-1}\sum_{m=0}^{w_2-1}\int_{y_0 = 0}^{t^*}f_{Y_0}(y_0)\frac{(t^*-y_0)^{n+m}}{n!m!}e^{-2t^*+2y_0}{\rm d}y_0 + \sum_{n=0}^{w_0-1}\frac{{t^*}^n}{n!}e^{-t^*}.
  \end{align*}
 Expanding $(t^*-y_0)^{n+m}$ yields
  \begin{align}
  &P(R(\bfa,j)>t^*) \nonumber \\
   &=\sum_{n=0}^{w_1-1}\sum_{m=0}^{w_2-1}\int_{y_0 = 0}^{t^*}\frac{y_0^{w_0-1}}{(w_0-1)!}e^{-y_0}\frac{1}{n!m!}\sum_{l=0}^{n+m}\binom{n+m}{l}{t^*}^l(-y_0)^{n+m-l}e^{-2t^*+2y_0}{\rm d}y_0 + \sum_{n=0}^{w_0-1}\frac{{t^*}^n}{n!}e^{-t^*} \nonumber\\
  &=\sum_{n=0}^{w_1-1}\sum_{m=0}^{w_2-1}\sum_{l=0}^{n+m}\frac{e^{-2t^*}{t^*}^l(-1)^{n+m-l}}{n!m!(w_0-1)!}\binom{n+m}{l}\int_{y_0 = 0}^{t^*}y_0^{n+m-l+w_0-1}e^{y_0}{\rm d}y_0 + \sum_{n=0}^{w_0-1}\frac{{t^*}^n}{n!}e^{-t^*}, \nonumber
\end{align}
completing the proof.
\end{proof}

\section{Additional numerical results}\label{sec:add_num}

\begin{table}
  \caption{Aggregate performance evaluated over 150 random graphs ($\rho = 0.04$, $\gamma = 0.6$)}
  \label{tbl:aggregate_results_expand2}
	\centering
  \includegraphics[width=\linewidth, trim={2.0cm 2.0cm 2.0cm 2.0cm}, clip]{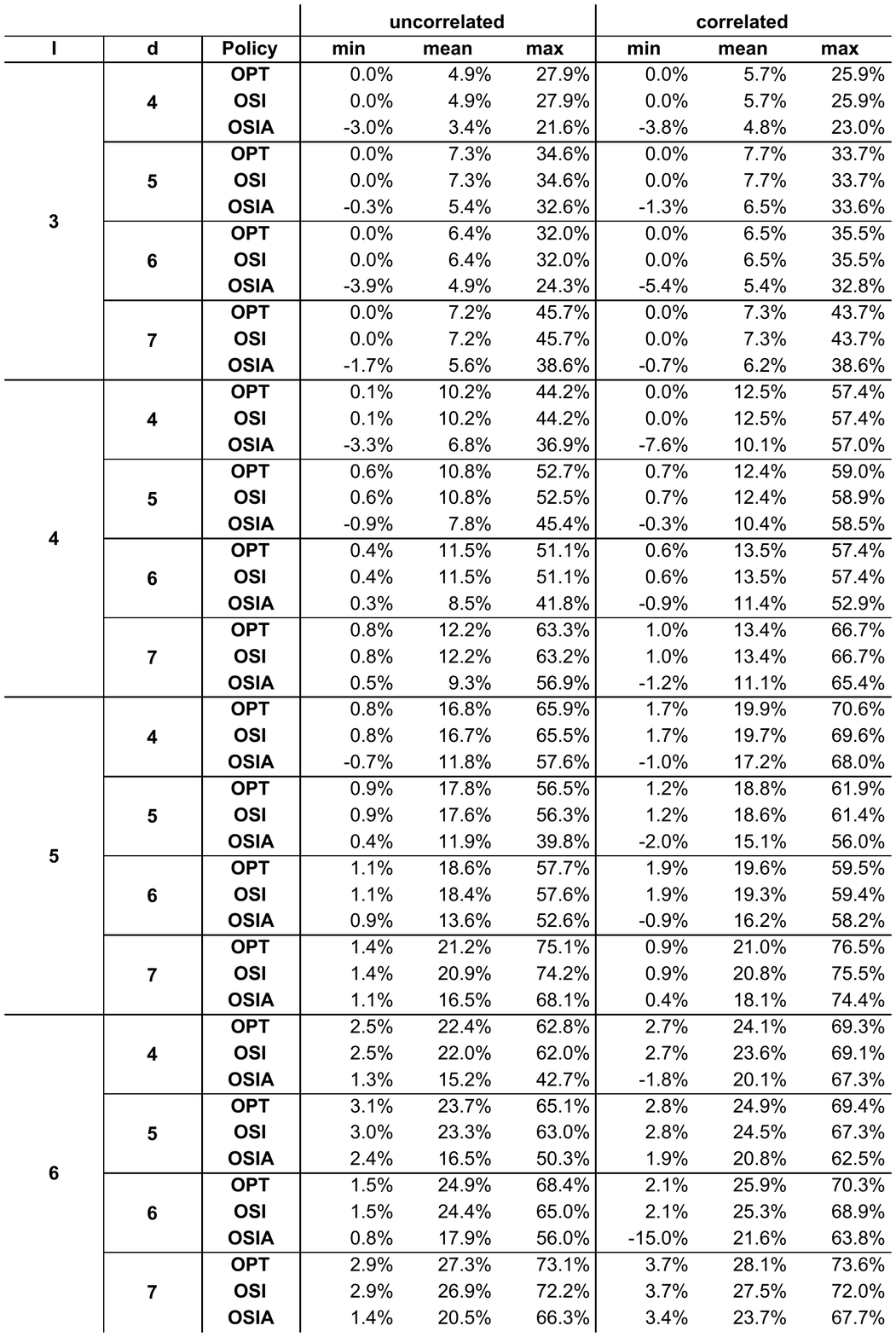}
\end{table}

\begin{table}
  \caption{Aggregate performance evaluated over 150 random graphs ($\rho = 0.02$, $\gamma = 0.6$)}
  \label{tbl:aggregate_results_expand3}
	\centering
  \includegraphics[width=\linewidth, trim={2.0cm 2.0cm 2.0cm 2.0cm}, clip]{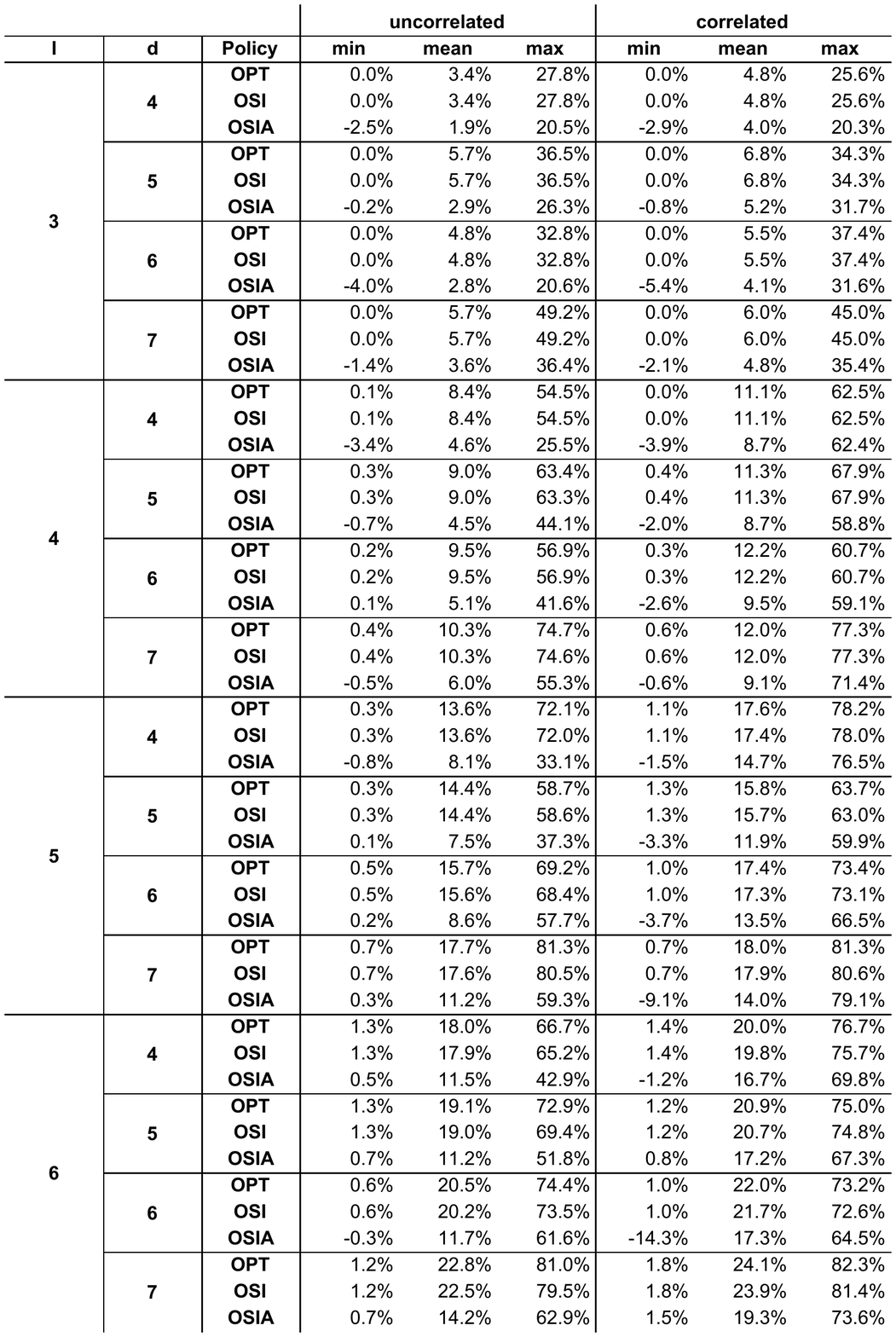}
\end{table}

\begin{table}
  \caption{Aggregate performance evaluated over multiple random graphs ($\rho = 0.1$, $\gamma = 0.6$)}
  \label{tbl:aggregate_results_expand5}
	\centering
  \includegraphics[width=\linewidth, trim={2.0cm 2.0cm 2.0cm 2.0cm}, clip]{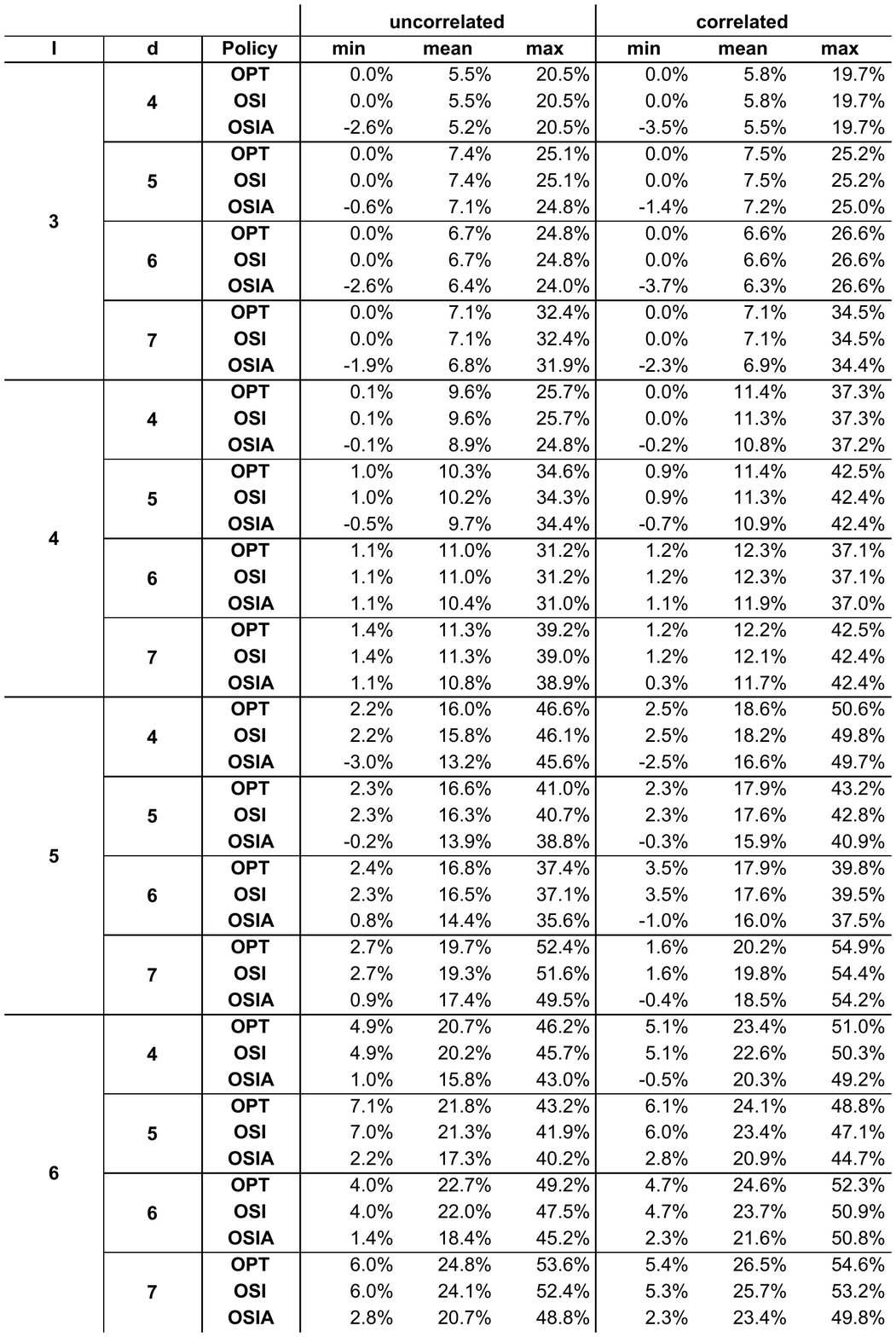}
\end{table}

\begin{table}
  \caption{Aggregate performance of OSIA over 50 random graphs ($\rho = 0.02$, $\gamma = 0.6$, $I=7$)}
  \label{tbl:aggregate_results_expand4}
	\centering
  \includegraphics[width=0.7\linewidth, trim={2.0cm 2.5cm 1.0cm 2.0cm}, clip]{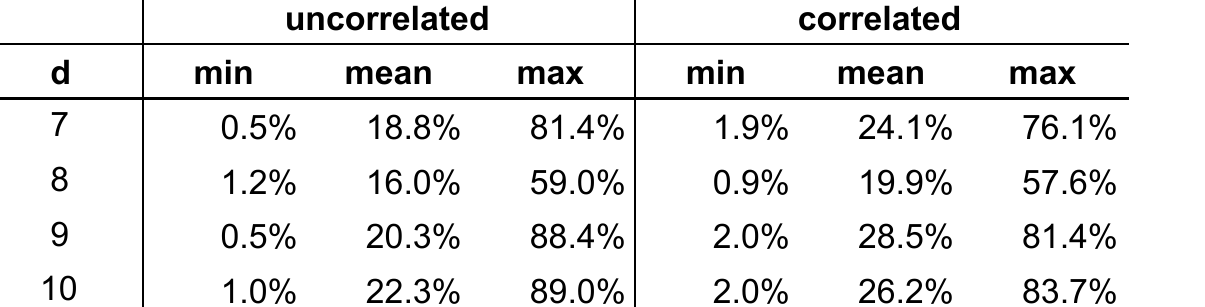}
\end{table}


\end{document}